\documentclass[a4paper,12pt, oneside]{scrartcl}

\usepackage{setspace}
\usepackage{tikz}
\usepackage{adjustbox}
\usepackage{tikz-cd}
\usetikzlibrary{arrows}
\usepackage{amssymb,amsmath,amsthm, dsfont}
\usepackage{bm}
\usepackage{mathtools}
\usepackage{upgreek}
\usepackage{graphicx}
\usepackage{amsfonts}
\usepackage{latexsym}
\usepackage[english]{babel}
\usepackage{verbatim} 
\usepackage[dvipsnames]{xcolor}
\usepackage{subcaption}
\usepackage{caption}
\usepackage{amsmath}
\usepackage{amsthm}
\usepackage{mathrsfs}
\usepackage{stmaryrd}
\usepackage{ragged2e}
\usepackage[bookmarks, colorlinks, breaklinks]{hyperref}
\hypersetup{
	linkcolor=black,
	citecolor=teal,
	filecolor=olive,
	urlcolor=olive }
\usepackage[a4paper,top=2.7cm,bottom=2.7cm,left=2.7cm,right=2.7cm]{geometry}
\usepackage{chapterbib}
\usepackage[T1]{fontenc}
\usepackage[utf8]{inputenc}
\usepackage{babel}
\usepackage{blindtext}

\newtheorem{lemma}{Lemma}[section]
\newtheorem{theorem}[lemma]{Theorem}
\newtheorem{proposition}[lemma]{Proposition}

\newtheorem*{theorem*}{Theorem}

\theoremstyle{definition}
\newtheorem{definition}[lemma]{Definition}

\newtheorem{notation}[lemma]{Notation}

\theoremstyle{remark}

\newtheorem{remark}[lemma]{Remark}

\newtheorem*{example*}{Example}

\newcommand{\Z}{\mathbb Z}
\newcommand{\Q}{\mathbb Q}
\newcommand{\Qb}{\overline{\mathbb Q}}

\newcommand{\Zp}{\mathbb Z_p}
\newcommand{\Fp}{\mathbb F_p}
\newcommand{\Fq}{\mathbb F_q}
\newcommand{\Qp}{\mathbb Q_p}

\newcommand{\Gm}{\mathbb G{_m}}

\newcommand{\AI}{\mathcal{A}_0}
\newcommand{\ApI}{\mathcal{A}_1}

\newcommand{\tr}{\text {Tr}}
\newcommand{\fr}{\text {Fr}}

\newcommand\blfootnote[1]{%
  \begingroup
  \renewcommand\thefootnote{}\footnote{#1}%
  \addtocounter{footnote}{-1}%
  \endgroup
}

\setkomafont{disposition}{\normalfont\bfseries}

\title{Trace of Frobenius on nearby cycles for a $\Gamma_1(p)$-local model}
\author{Giulio Marazza}
\date{}

\begin{document}

\hypersetup{
	linkcolor=Mahogany,
	citecolor=teal,
	filecolor=olive,
	urlcolor=olive }

\maketitle

\begin{abstract}
    \small{\textsc{ABSTRACT}}.\, We compute the trace of Frobenius on the sheaf of nearby cycles for the integral model of the Siegel modular variety with pro-$p$ Iwahori level structure, constructed by Haines and Stroh, in the case of $\text{GSp}_4$. To this end, we make use of a semistable resolution of this model due to Shadrach. Finally, we compare our results with calculations made by Horn, confirming a conjecture of Haines and Kottwitz in this case.
\end{abstract}

\section{Introduction}

Within the context of the Langlands program, it is an important task to determine the local factor at a prime $p$ of the Hasse-Weil zeta function of a given Shimura variety, with the hope of expressing this function as a product of automorphic $L$-functions. In the case of good reduction, inertia acts trivially on the étale cohomology of the Shimura variety and the smooth base change theorem holds. The situation is more complicated in the case of bad reduction, as these facts are no longer true. To take into account the inertia action, one is led to work with the semisimple zeta function, for which we refer to \cite{HN}; to remedy the failure of the smooth base change theorem we need to consider the étale cohomology of the sheaf of nearby cycles on the special fiber of a suitable integral model of the Shimura variety. More precisely, one is led to consider the alternating sum of the semisimple trace  of Frobenius on the stalks of nearby cycles at points of the special fiber.\\ In this paper, we are concerned with the Siegel modular variety with pro-$p$ Iwahori level structure in the case of $\text{GSp}_4$, where explicit calculations are still feasible. We consider the integral model constructed by Haines and Stroh in \cite{HS}, denoted by $\mathcal{A}_1$: it is defined over the integral model with Iwahori level structure, written $\mathcal{A}_0$, which is the moduli space of chains of degree $p$ isogenies between abelian varieties of dimension two. The model $\mathcal{A}_1$ then parametrizes Oort-Tate generators (see Definition \ref{OT generators}) for the kernel of these isogenies. We compute the value of 
\[
\text{Tr}^\text{ss}(\text{Fr}_q; (\pi_*R\Psi^1)_x)=\sum_i(-1)^i\text{Tr}^\text{ss}(\text{Fr}_q; (\pi_*R^i\Psi^1)_x),
\]
\blfootnote{giulio.marazza@uni-due.de}
where $\pi:\mathcal{A}_1\rightarrow \mathcal{A}_0$ is the natural morphism from the $\Gamma_1(p)$-model to the $\Gamma_0(p)$-model, $R\Psi^1$ is the complex of nearby cycles on $\mathcal{A}_{1,\overline{\mathbb{F}}_p}$ and $q$ is a power of $p$. We then compare it with the function $\phi_{r,1}'$ constructed by Horn in \cite{Hor} for $\text{GSp}_4$, proving the following result.

\begin{theorem*}[\ref{comparison test function}]
    For all $x\in \mathcal{A}_0(\Fq)$ we have that
    \[
    \emph{Tr}^\emph{ss}(\emph{Fr}_q; (\pi_*R\Psi^1)_x)=(q-1)^3\phi'_{r,1}(s_xw),
    \]
    where $w\in \emph{Adm}(\mu)$ corresponds to the KR-stratum which $x$ belongs to. 
\end{theorem*}

Here, $s_x$ is an element of $T(\Fq)$, where $T\subset \text{GSp}_4$ is the standard maximal torus, which depends on the point $x$ and some extra data. The function $\phi_{r,1}'$ is an element of the pro-$p$ Iwahori-Hecke algebra for $\text{GSp}_4$ , so it is determined by its values on the elements of $\widetilde{W}_1=T(\Fq)\rtimes \widetilde{W}$, where $\widetilde{W}$ is the extended affine Weyl group. Finally, $\text{Adm}(\mu)\subset \widetilde{W}$ is the $\mu$-admissible subset for the cocharacter $\mu=(1,1,0,0)$ of $\text{GSp}_4$, whose elements parametrize the Kottwitz-Rapoport stratification on $\mathcal{A}_{0,\Fp}$. The function $\phi'_{r,1}$ is a modification of another element of the pro-$p$ Iwahori-Hecke algebra, $\phi_{r,1}$, defined in \cite{Hai}, which should conjecturally enter the counting point formula needed to relate the Hasse-Weil zeta function of the Shimura variety we are considering with automorphic $L$-functions. \\
We obtain our result by intricate, though elementary, computations on the local model of $\ApI$, denoted by ${U}_1$: this is a scheme connected to $\mathcal{A}_1$ by a smooth correspondence, therefore the nearby cycles on $\mathcal{A}_1$ and ${U}_1$ have isomorphic stalks on corresponding points. The local model that we are going to use is constructed in \cite{HS}, see also \cite{Mar} for a more conceptual description of it and for some of its geometric properties. The advantage is that $U_1$ can be explicitly described in terms of generators and equations. The computation cannot be directly performed on $U_1$, due to the presence of singularities: we have to consider a semistable resolution due to Shadrach in \cite[Chapter 6]{Sha}, consisting of a sequence of blow-ups with centers in the special fiber $U_{1,\Fp}$, which we make explicit. \\
This kind of result is also proved by Haines and Rapoort in \cite{HR}, where they consider the Drinfeld case for unitary Shimura varieties and prove stronger statements. The situation here presents similarities with their case, but also some differences. Most notably, the element $s_x$ in their setting can be always determined in terms of the Oort-Tate parameters of the point $x$ (see Section \ref{comparison test function}, in particular Remark \ref{specific element}), while in the $\text{GSp}_4$ case this does not hold for a specific $KR$-stratum, see the discussion about $w=s_{02}\tau$ in Subsection \ref{the GSp4 case}. It would be interesting to find a more conceptual description of this element $s_x$. Let us also mention ongoing work of Haines, Li, Stroh \cite{HLS} and of Li \cite{Li} towards the proof of this conjecture, at least in the linear and symplectic case. \\
Here is the layout of the paper. In Chapter \ref{The integral models} we introduce the relevant integral models with their associated local model diagrams, both in the Iwahori and pro-$p$ Iwahori case. In Chapter \ref{A semistable resolution of the local models} we describe the semistable resolutions of both local models, the one in the Iwahori case due to de Jong in \cite{dJo}. In Chapter \ref{nearby cycles} we compute the trace of Frobenius on nearby cycles for $U_1$. In Chapter \ref{Comparison with the test function} we recall the work of Haines and Rapoport in the good Drinfeld case and then we compare our calculations with the ones made by Horn.

\bigskip
\bigskip

I would like to thank Ulrich Görtz for his constant interest and encouragement in my work. This manuscript is part of my PhD thesis, which was carried out in the DFG Research Training Group 2553. I would also like to thank Qihang Li for a helpful comment.

\section{The integral models}\label{The integral models}

In this section, we recall the construction of the integral models of the Siegel modular variety for $\text{GSp}_{2n}$ with $\Gamma_0(p)$ and $\Gamma_1(p)$-level structure that we are going to consider, together with their local models. We will specialize to the $\text{GSp}_4$ case in the next section.

\begin{definition}[\cite{RZ}, Definition 6.9]\label{GSp iwahori integral}
    The integral model for the Siegel modular variety with Iwahori level structure at $p$ is the scheme $\mathcal{A}_0$ representing the functor which sends a $\Z_p$-scheme $S$ to the set of isomorphism classes of tuples $(A_\bullet, \lambda_0, \lambda_n, \overline{\eta})$, where:
    \begin{itemize}
        \item [$(i)$]$A_\bullet$ is a chain $A_0\xrightarrow{\alpha_0}A_1\xrightarrow{\alpha_1}\cdots\xrightarrow{\alpha _{n-1}}A_n$ of $n$-dimensional abelian schemes over $S$, with $\alpha_i$ isogenies of degree $p$;
        \item [$(ii)$] the maps $\lambda_0, \lambda_n$ are principal polarizations of $A_0, A_n$ respectively, making the loop starting at any $A_i$ or $A_i^\vee$ in the diagram
        \[
        \begin{tikzcd}
            A_0 \arrow[r, "\alpha_0"] & A_1 \arrow[r,"\alpha_1"] & \cdots \arrow[r,"\alpha_{n-1}"] & A_{n} \arrow[d,"\lambda_n"] \\
            A_0^\vee \arrow[u, "\lambda_0^{-1}"] & A_1^\vee \arrow[l, "\alpha_0^\vee"] & \cdots \arrow[l, "\alpha_1^\vee"] & A_{n}^\vee \arrow[l, "\alpha_{n-1}^\vee"]
        \end{tikzcd}
        \]
        multiplication by $p$;
        \item [$(iii)$] $\overline{\eta}$ is a $K^p$-level structure.
    \end{itemize}
\end{definition}

\begin{notation}
    Given an $S$-point $(A_\bullet, \lambda_0, \lambda_n, \overline{\eta})$ of $\mathcal{A}_0$, we are going to write $G_i$ for the kernel of $\alpha_i$. 
\end{notation}

Let now $V=\Qp^{2n}$, $\{e_1,\ldots, e_{2n}\}$ the canonical basis and $(\cdot,\cdot)$ the standard symplectic pairing given by $(e_i,e_{2n+1-j})=\delta_{ij}$. Consider the full self-dual $\Zp$-lattice chain $\Lambda_\bullet$ given by
\[
\Lambda_i=\langle p^{-1}e_1,\ldots, p^{-1}e_i, e_{i+1},\ldots, e_{2n} \rangle_{\Zp}
\]
for $i=0,\ldots, 2n-1$ and extended periodically for all other $i\in\Z$. With this notation, the local model for $\mathcal{A}_0$ is the scheme $\text{M}_0$ whose $S$-valued points, for a $\Zp$-scheme $S$, is the set of commutative diagrams
\[
\begin{tikzcd}
{\Lambda_{0,S}} \arrow[r] & {\Lambda_{1,S}} \arrow[r] & \cdots \arrow[r] & {\Lambda_{2n-1,S}} \arrow[r] & p^{-1}\Lambda_{0,S} \\
F_0 \arrow[u, hook] \arrow[r] & F_1 \arrow[u, hook] \arrow[r] & \cdots \arrow[r] & F_{2n-1} \arrow[u, hook] \arrow[r] & p^{-1}F_0 \arrow[u, hook]
\end{tikzcd}
\]
where each $F_i$ is a locally free $\mathcal{O}_S$-module of rank $n$ which is Zariski-locally a direct summand of $\Lambda_{i,S}$. Moreover, for each $i\in\{ 0,\ldots,2n-1 \}$, the map 

\begin{equation}\label{condition loc mod iw}
F_i\hookrightarrow \Lambda_{i,S}\cong\Lambda^\vee_{2n-i,S}\rightarrow F_{2n-i}^\vee
\end{equation}

is required to be an isomorphism. The integral model $\mathcal{A}_0$ and the local model $\text{M}_0$ are connected by a smooth correspondence of the form

\begin{equation}\label{loc model diagram Iw}
\begin{tikzcd}
    & \tilde{\mathcal{A}_0} \arrow[ld, swap, "\phi"] \arrow[rd, "\psi"] & \\
    \text{M}_0 & & {\mathcal{A}_0},
\end{tikzcd}
\end{equation}

where $\tilde{\mathcal{A}_0}$ parametrizes points of $\mathcal{A}_0$ together with an isomorphism $\gamma_\bullet$ between the chain of the de Rham cohomology groups $H^1_{dR}(A_\bullet)$ and $\Lambda_\bullet$. The automorphism group scheme $\mathcal{G}=\underline{\text{Aut}}(\Lambda_\bullet)$ acts on $\tilde{A}_0$, making $\psi$ a $\mathcal{G}$-torsor and $\phi$ a $\mathcal{G}$-equivariant map. As in \cite[Chapter 5]{Goer2}, we can embed the special fiber of $\text{M}_0$ into the affine flag variety for $\text{GSp}_{2n}$, realizing $\text{M}_{0,\Fp}$ as a union of Schubert cells:
\[
\text{M}_{0,\Fp}=\coprod_{w\in\text{Adm}(\mu)} S_w.
\]
This is the Kottwitz-Rapoport (or KR) stratification. Here $\text{Adm}(\mu)$ is the $\mu$-admissible set, where $\mu$ is the cocharacter of $\text{GSp}_{2n}$ given by $\mu=(1^{(n)},0^{(n)})$. It is a subset of the extended affine Weyl group $\widetilde{W}$ for $\text{GSp}_{2n}$, defined as
\[
\text{Adm}(\mu)=\{ w\in \widetilde{W} \;|\; \exists\, \tau \in W \; \text{such that} \; w\leq t_{\tau(\mu)} \};
\]
here, for an element $w\in\widetilde{W}$, $t_w$ stands for its translation part. Each $S_w$ is a $\mathcal{G}$-orbit, so we can descend this stratification to $\mathcal{A}_{0,\Fp}$, and we are going to write $\mathcal{A}_{0,w}$ for the strata on the integral model side. By the main result of \cite{KR} in the symplectic case, one can identify $\text{Adm}(\mu)$ with the set of permissible alcoves, defined as follow.

\begin{definition}[\cite{KR}, Section 4.2]
    An alcove $x=(x_i)_i$ for $\text{GSp}_{2n}$ is a collection of vectors $x_0, \cdots, x_{2n-1}\in \Z^{2n}$ such that:
    \begin{itemize}
        \item[$(i)$] $x_0\leq x_1\leq \dots \leq x_{2n-1}\leq x_n=x_0 + \textbf{1}$, where inequality is component-wise and $\textbf{1}=(1, \ldots, 1)\in \Z^{2n}$;
        \item[$(ii)$] $\sum_j x_{i+1}(j) = \sum_j x_{i}(j) +1$ for each $i=0,\ldots,2n-1$;
        \item [$(iii)$] there exists $d\in\Z$ such that $x_{2n-1-i}=\textbf{d} + \theta(x_i)$ for $i=0,\ldots, 2n-1$. Here $\textbf{d}=(d,\ldots,d)\in \Z^{2n}$ and $\theta$ is the automorphism of $\mathbb{R}^{2n}$ defined by
    \[
    \theta(x_1,x_2,\ldots, x_{2n})=(-x_{2n},\ldots, -x_2, -x_1).
    \]
    \end{itemize}
    An alcove $x=(x_i)_i$ is called permissible if it satisfies the following two conditions:
    \begin{itemize}
        \item [$(a)$] $\sum_j x_0(j)=n$;
        \item [$(b)$] $\omega_i \leq x_i \leq \omega_i +\textbf{1} $, where $\omega_i=(1^{(i)}, 0^{(2n-i)})\in \Z^{2n}$, for all $i=0,\ldots, 2n-1$.
    \end{itemize}
    We define its difference vectors to be $t_i^x=x_i-\omega_i$.
\end{definition}

For each permissible alcove $x$, following \cite[Definition 4.4]{Goer}, we define the open subscheme $U_x\subset\text{M}_0$ whose functor of points is
\begin{equation}\label{open unitary}
    U_x(R)=\Bigg\{ F_\bullet \in \text{M}_0(R) \;|\; F_i\oplus \Big( \bigoplus_{t_i^x(j)=0} R\text{e}_j\Big)=R^n \Bigg\}:
\end{equation}
its special fiber contains the stratum $S_x$. Moreover, there exists a unique element $\uptau\in \text{Adm}(\mu)$ which is minimal with respect to the Bruhat-Tits length on $\widetilde{W}$ and whose corresponding stratum is a single point, called the worst singular point of the local model and written by $\tau$. The $\mathcal{G}$-translates of $U_\uptau$ cover the whole $\text{M}_0$: this means that we can use $U_\uptau$ itself as a local model for $\mathcal{A}_0$. For the rest of this text, we are going to write $U$ instead of $U_\uptau$. It is an affine scheme, whose coordinate ring $B$ admits an explicit presentation, see \cite[Chapter 4]{Sha}.\\ In order to construct the integral model $\mathcal{A}_1$, we need the classification of finite locally free group schemes of rank $p$ over a $\Zp$-base due to Tate and Oort from \cite{OT}. We recall a summary of their main results as in \cite[Theorem 3.3.1]{HR}.

\begin{theorem}[Tate-Oort]
    Let $OT$ be the $\Zp$-stack of finite locally free group schemes of rank $p$.
    \begin{itemize}
        \item [$(i)$] OT is an Artin stack isomorphic to
        \[
            OT=\left[ \emph{Spec}\left(\frac{\Zp[x,y]}{(xy-\omega_p)}\right)/\Gm \right]
        \]
        where $\Gm$ acts by $t\cdot (x,y)=(t^{p-1}x, t^{1-p}y)$ and $\omega_p$ is an explicit element of $p\Z_p^\times$.
        
        \item [$(ii)$] The universal group scheme $\mathcal{G}_{OT}$ over OT is given by
        \[
            \mathcal{G}_{OT}=\left[ \emph{Spec}\left(\frac{\Zp[x,y,z]}{(xy-\omega_p, z^p-xz)}\right)/\Gm \right]
        \]
        where $\Gm$ acts as before on $x, y$ and with weight $1$ on $z$, with zero section $z=0$.

        \item[$(iii)$] Cartier duality acts on $OT$ by interchanging $x$ and $y$.
    \end{itemize}
\end{theorem}

Now consider the closed substack
\[
\mathcal{G}_{OT}^\times = \left[ \text{Spec}\left(\frac{\Zp[x,y,z]}{(xy-\omega_p, z^{p-1}-x)}\right)/\Gm \right]
\]
of $\mathcal{G}_{OT}$, called the substack of generators of $\mathcal{G}_{OT}$.

\begin{definition}[\cite{HR}, Section 3.3]\label{OT generators}
    Let $S$ be a $\Z_p$-scheme, $G$ a finite locally free group scheme of rank $p$ over $S$. Let $G^\times$ be the closed subscheme of $G$ obtained by pulling back $\mathcal{G}_{OT}^\times\subset \mathcal{G}_{OT}$ along the map $S\xrightarrow{}OT$ corresponding to $G$. A global section $z\in G(S)$ is called an {Oort-Tate generator} if it factors through $G^\times$.
\end{definition}

\begin{remark}
    Note that a map $S\rightarrow OT$ is the same as a triple $(L,a,b)$ where $L$ is a line bundle on $S$, $a\in\Gamma(S,L^{p-1})$ and $b\in \Gamma(S,L^{1-p})$ are sections such that $a\otimes b=\omega_p$: then, an Oort-Tate generator is a section $x\in\Gamma(S,L)$ such that $x^{\otimes p-1}=a$. Given $G$ a finite locally free $S$-group scheme of rank $p$ and $g, g^*$ Oort-Tate generators for $G$ and its Cartier dual $G^*$, the product $g^*\otimes g$ is a global section of $\mathcal{O}_S$ since the line bundle corresponding to $G^*$ is the dual of the one corresponding to $G$. The line bundle corresponding to $G$ will be called the Oort-Tate line bundle.
\end{remark}

\begin{definition}[\cite{HS}]
    The integral model for the Siegel modular variety with pro-$p$ Iwahori level structure at $p$ is the scheme $\mathcal{A}_1$ representing the functor which sends a $\Z_p$-scheme $S$ to the set of isomorphism classes of tuples $\big( (A_\bullet, \lambda_0, \lambda_n, \overline{\eta}), (g_i, g^*_i)_{i=0}^{n-1} \big)$, where:
    \begin{itemize}
        \item [$(i)$] $(A_\bullet, \lambda_0, \lambda_n, \overline{\eta})$ is an $S$-point of $\mathcal{A}_0$;
        \item [$(ii)$] $g_i$ and $g_i^*$ are Oort-Tate generators of $G_i$ and $G_i^*$ respectively, such that $g_i^*g_i$ is independent of $i$.
    \end{itemize}
\end{definition}

We choose as a $\Gamma_1(p)$-local model the affine scheme $U_1=\text{Spec}(C)$, where
\[
C= \frac{B[u_0,v_0,\ldots,u_{n-1},v_{n-1}]}{\left(\substack{u_i^{p-1}-\alpha_i,\, v_i^{p-1}-\beta_i\\ u_0v_0-u_iv_i}\right)_i}.
\]
Here, $\alpha_i$ and $\beta_i$ are suitable elements of $B$ which correspond, under the local model diagram (\ref{loc model diagram Iw}), to the Oort-Tate parameters of the group schemes $G_i$: in the case of $\text{GSp}_4$ we are going to make them explicit. See \cite{Mar} for more details about how the local model diagram works in the $\Gamma_1(p)$-setting.

\section{A semistable resolution of the local models\label{A semistable resolution of the local models}}

We review the resolution of the Iwahori local model $U$ for $\text{GSp}_4$ constructed by de Jong in \cite{dJo}, since it will be the first step for the resolution of $U_1$. We also compute the equations of each KR-stratum, since we are going to need them later.

\subsection{Iwahori case}

In the $\text{GSp}_4$ case, the alcove corresponding to the worst singular point has difference vectors $(0011),(1001),(1100),(0110)$. The $R$-points of $U$ are tuples $(F_0,F_1,F_2,F_3)$ sitting in the following diagram
\[
\begin{tikzcd}[row sep=1em]
R^4 \arrow[r, "\varphi_0"]    & R^4 \arrow[r, "\varphi_1"]    & R^4 \arrow[r, "\varphi_2"]    & R^4 \arrow[r, "\varphi_3"]    & R^4                 \\
F_0 \arrow[r] \arrow[u, hook] & F_1 \arrow[r] \arrow[u, hook] & F_2 \arrow[r] \arrow[u, hook] & F_3 \arrow[r] \arrow[u, hook] & F_0, \arrow[u, hook]
\end{tikzcd}
\]
where $\varphi_i=\text{diag}(1,\dots, p,\dots, 1)$ with $p$ in the $(i+1)$-th component. They are subject to the condition (\ref{condition loc mod iw}), from which it follows that $F_1$ determines $F_3$ and that $F_0,F_2$ are isotropic with respect to $(\cdot,\cdot)$, and to the condition (\ref{open unitary}), from which it follows that the $F_i$ are free $R$-modules of rank 2. Choosing suitable basis for $F_0,F_1$ and $F_2$ we can express the $R$-points of $U$ as chains of the form

\begin{equation}\label{matr equ g=2}
        \begin{tikzcd} 
        \setlength\arraycolsep{1pt}
        \scalebox{0.9}{ $\begin{pmatrix}
            1 & 0 \\   
            0 & 1 \\
            a^0_{11}  & a^0_{12}\\
            a^0_{21}  & a^0_{11}
        \end{pmatrix} $} \arrow[r, "\varphi_{0|F_0}"] 
        &
        \setlength\arraycolsep{1pt}
        \scalebox{0.9}{ $\begin{pmatrix}
            a^1_{21}  & a^1_{22}\\
            1 & 0 \\   
            0 & 1 \\
            a^1_{11}  & a^1_{12}
        \end{pmatrix} $}  \arrow[r, "\varphi_{1|F_1}"]
        &
        \setlength\arraycolsep{1pt}
        \scalebox{0.9}{ $\begin{pmatrix}
            a^2_{11}  & a^2_{12}\\
            a^2_{21}  & a^2_{11}\\
            1 & 0 \\   
            0 & 1            
        \end{pmatrix} $} 
    \end{tikzcd}
\end{equation}

and imposing the condition $\varphi_i(F_i)\subset F_{i+1}$ gives the relations between the generators $a^i_{jk}$. It is easy to check that all the variables can be expressed in terms of $a^0_{11}, a^1_{22}, a^2_{12}, a^0_{12}, a^1_{12}$, obtaining a presentation of $U$ as the spectrum of

\begin{equation}\label{loc mod g=2}
B=\Zp[x,y,a,b,c]/(xy-p, ax+by+abc).
\end{equation}

Blowing up $U$ along the closed subscheme given by $(x,b)$ yields a semistable resolution, which can be described as
\[
U'=\text{Proj}_{U}\big( B[\Tilde{x},\Tilde{b}]/(x\Tilde{b}-b\Tilde{x}, a\Tilde{x}+\Tilde{b}y+a\Tilde{b}c) \big),
\]
with the two standard charts given by:

\begin{itemize}
    \itemsep-1em
    \item[]
    \[
    U'_{[\tilde{b}=1]}=\text{Spec}\left( \frac{\mathbb{Z}[\tilde{x}, a, b, c]}{(\tilde{x}abc - p)} \right)
    \longrightarrow \text{Spec}(B),
    \quad
    \left.
    \begin{array}{l}
    b\tilde{x} \\[-5pt]
    ac \\[-5pt]
    -a \\[-5pt]
    c - \tilde{x}
    \end{array}
    \right\}
    \longmapsfrom
    \left\{
    \begin{array}{l}
    x \\[-4pt]
    y \\[-4pt]
    a \\[-4pt]
    c
    \end{array}
    \right. ;
    \]
    
    \item[]
    \[
    U'_{[\tilde{x}=1]}=\text{Spec}\left( \frac{\mathbb{Z}[x,y, a, \tilde{b}, c]}{(xy - p, a+\tilde{b}y+a\tilde{b}c)} \right)
    \longrightarrow \text{Spec}(B), \quad \tilde{b}x \longmapsfrom b;
    \]
\end{itemize}

The chart $U'_{[\tilde{x}=1]}$ is covered by the two opens where $\tilde{b}$, $1+\tilde{b}c$ are invertible, over which it takes the following form:

\begin{itemize}
        \itemsep-1em
        \item[]
        \[
        U'_{[\tilde{x}=1, (1+\tilde{b}c)^{-1}]}=\text{Spec}\left( \frac{\mathbb{Z}[x,y, \tilde{b}, c,(1+\tilde{b}c)^{-1}]}{(xy - p)} \right)\longrightarrow \text{Spec}(B), \quad \frac{-\tilde{b}y}{1+\tilde{b}c}\longmapsfrom a;
        \]

        \item[]
        \[
        U'_{[\tilde{x}=1, \tilde{b}^{-1}]}=\text{Spec}\left( \frac{\mathbb{Z}[x,a, \tilde{b}^\pm, c]}{(xac - p)} \right)\longrightarrow \text{Spec}(B),
        \quad
        \left.
        \begin{array}{l}
        x \\[-7pt]
        ac \\[-7pt]
        a\\[-7pt]
        -\tilde{b}x \\[-7pt]
        -c+ \tilde{b}^{-1}
        \end{array}
        \right\}
        \longmapsfrom
        \left\{
        \begin{array}{l}
        x \\[-4pt]
        y \\[-4pt]
        a \\[-4pt]
        b \\[-4pt]
        c
        \end{array}
        \right. .
        \]
\end{itemize}

The KR-stratification of $U_{ \Fp}$ has 13 strata $S_w$, where $w\in \text{Adm}(\mu)$ for $\mu=(1,1,0,0)$, and $U_{\Fp}$ has four irreducible components, obtained as the closures of the four open strata indexed by the translation elements $t_\lambda\in \text{Adm}(\mu)$ (see the picture below). We are going to need an explicit description of each stratum in terms of the equations (\ref{loc mod g=2}), which we can get by working with the opens $V_w=U\cap U_w$ for $w\in \text{Adm}(\mu)$: if $\ell(w)=3$, $S_w=V_{w,\Fp}$; if $\ell(w)=2$, $S_w=(V_{w,\Fp})^\text{sing}$; if $\ell(w)=1$ we need to intersect the irreducible components containing the stratum and remove the worst singular point, which is given by $\uptau=V(x,y,a,b,c)$. From the definition of the opens $U_w$ in (\ref{open unitary}), we see that a point $F_\bullet\in U$ lies also in $U_w$ if and only if the condition
\[
F_i\oplus \Big( \bigoplus_{t_i^w(j)=0} R\text{e}_j\Big)=R^n
\]
is satisfied, where $t_i^w$ are the difference vectors of the alcove corresponding to $w\in\text{Adm}(\mu)$. But this happens if and only if the minor of $F_i$ given by the rows where the coordinates of $t_i^w$ are 1 is invertible.

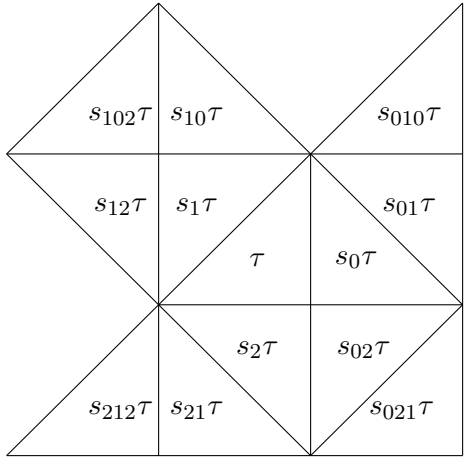
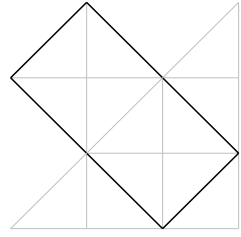
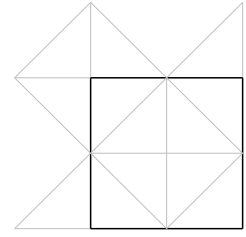
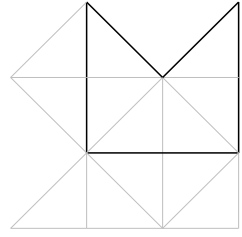
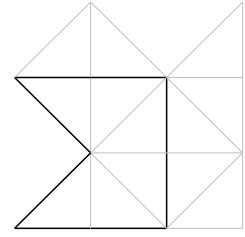
\adjustbox{scale=1, center}{
\begin{tikzpicture}

\node at (3.3,2.6) {\small $\tau$};
\node at (4.6,2.6) {\small $s_0\tau$};
\node at (2.5,3.3) {\small $s_1\tau$};
\node at (3.3, 1.4) {\small $s_2\tau$};
\node at (4.7, 1.4) {\small $s_{02}\tau$};
\node at (1.5,3.3) {\small $s_{12}\tau$};
\node at (1.5,4.5) {\small $s_{102}\tau$};
\node at (2.5,4.5) {\small $s_{10}\tau$};
\node at (5.3,3.3) {\small $s_{01}\tau$};
\node at (5.3,4.5) {\small $s_{010}\tau$};
\node at (5.2, 0.6) {\small $s_{021}\tau$};
\node at (2.5, 0.6) {\small $s_{21}\tau$};
\node at (1.5, 0.6) {\small $s_{212}\tau$};
\draw (0,0) -- (6,6);
\draw (0,0) -- (6,0);
\draw (6,0) -- (6,6);
\draw (4,0) -- (0,4);
\draw (2,6) -- (0,4);
\draw (2,6) -- (6,2);
\draw (2,6) -- (2,0);
\draw (0,4) -- (6,4);
\draw (4,0) -- (6,2);
\draw (4,0) -- (4,4);
\draw (2,2) -- (6,2);

\node at (9.7,-0.8) {\tiny{$t_{(1100)}, \, V(x,b)$}};
\draw[line width=0.2mm] (11,2.5) -- (11,1.5);
\draw[line width=0.2mm] (11,1.5) -- (11,0.5);
\draw[line width=0.11mm,lightgray] (8,-0.5) -- (9,-0.5);
\draw[line width=0.11mm,lightgray] (9,-0.5) -- (10,-0.5);
\draw[line width=0.11mm,lightgray] (10,-0.5) -- (11,-0.5);
\draw[line width=0.11mm,lightgray] (11,-0.5) -- (11,0.5);
\draw[line width=0.11mm,lightgray] (9,-0.5) -- (9,0.5);
\draw[line width=0.11mm,lightgray] (10,-0.5) -- (10,1.5);
\draw[line width=0.2mm] (9,0.5) -- (11,0.5);
\draw[line width=0.2mm] (9,1.5) -- (9,2.5);
\draw[line width=0.2mm] (9,0.5) -- (9,1.5);
\draw[line width=0.11mm,lightgray] (8,1.5) -- (9,1.5);
\draw[line width=0.11mm,lightgray] (9,1.5) -- (10,1.5);
\draw[line width=0.11mm,lightgray] (10,1.5) -- (11,1.5);
\draw[line width=0.11mm,lightgray] (8,1.5) -- (9,0.5);
\draw[line width=0.11mm,lightgray] (9,0.5) -- (10,-0.5);
\draw[line width=0.2mm] (9,2.5) -- (10,1.5);
\draw[line width=0.11mm,lightgray] (10,1.5) -- (11,0.5);
\draw[line width=0.11mm,lightgray] (8,-0.5) -- (9,0.5);
\draw[line width=0.11mm,lightgray] (8,1.5) -- (9,2.5);
\draw[line width=0.11mm,lightgray] (9,0.5) -- (10,1.5);
\draw[line width=0.11mm,lightgray] (10,-0.5) -- (11,0.5);
\draw[line width=0.2mm] (10,1.5) -- (11,2.5);

\node at (13.55,-0.8) {\tiny{$t_{(0011)}, \, V(y,a)$}};
\draw[line width=0.11mm,lightgray] (15,2.5) -- (15,1.5);
\draw[line width=0.11mm,lightgray] (15,1.5) -- (15,0.5);
\draw[line width=0.2mm] (12,-0.5) -- (13,-0.5);
\draw[line width=0.2mm] (13,-0.5) -- (14,-0.5);
\draw[line width=0.11mm,lightgray] (14,-0.5) -- (15,-0.5);
\draw[line width=0.11mm,lightgray] (15,-0.5) -- (15,0.5);
\draw[line width=0.11mm,lightgray] (13,-0.5) -- (13,0.5);
\draw[line width=0.2mm] (14,-0.5) -- (14,1.5);
\draw[line width=0.11mm,lightgray] (13,0.5) -- (15,0.5);
\draw[line width=0.11mm,lightgray] (13,1.5) -- (13,2.5);
\draw[line width=0.11mm,lightgray] (13,0.5) -- (13,1.5);
\draw[line width=0.2mm] (12,1.5) -- (13,1.5);
\draw[line width=0.2mm] (13,1.5) -- (14,1.5);
\draw[line width=0.11mm,lightgray] (14,1.5) -- (15,1.5);
\draw[line width=0.2mm] (12,1.5) -- (13,0.5);
\draw[line width=0.11mm,lightgray] (13,0.5) -- (14,-0.5);
\draw[line width=0.11mm,lightgray] (13,2.5) -- (14,1.5);
\draw[line width=0.11mm,lightgray] (14,1.5) -- (15,0.5);
\draw[line width=0.2mm] (12,-0.5) -- (13,0.5);
\draw[line width=0.11mm,lightgray] (12,1.5) -- (13,2.5);
\draw[line width=0.11mm,lightgray] (13,0.5) -- (14,1.5);
\draw[line width=0.11mm,lightgray] (14,-0.5) -- (15,0.5);
\draw[line width=0.11mm,lightgray] (14,1.5) -- (15,2.5);

\node at (13.7, 3.2) {\tiny{$t_{(1010)}, \, V(x, y+ac)$}};
\draw[line width=0.11mm,lightgray] (15,6.5) -- (15,5.5);
\draw[line width=0.2mm] (15,5.5) -- (15,4.5);
\draw[line width=0.11mm,lightgray] (12,3.5) -- (13,3.5);
\draw[line width=0.2mm] (13,3.5) -- (14,3.5);
\draw[line width=0.2mm] (14,3.5) -- (15,3.5);
\draw[line width=0.2mm] (15,3.5) -- (15,4.5);
\draw[line width=0.2mm] (13,3.5) -- (13,4.5);
\draw[line width=0.11mm,lightgray] (14,3.5) -- (14,5.5);
\draw[line width=0.11mm,lightgray] (13,4.5) -- (15,4.5);
\draw[line width=0.11mm,lightgray] (13,5.5) -- (13,6.5);
\draw[line width=0.2mm] (13,4.5) -- (13,5.5);
\draw[line width=0.11mm,lightgray] (12,5.5) -- (13,5.5);
\draw[line width=0.2mm] (13,5.5) -- (14,5.5);
\draw[line width=0.2mm] (14,5.5) -- (15,5.5);
\draw[line width=0.11mm,lightgray] (12,5.5) -- (13,4.5);
\draw[line width=0.11mm,lightgray] (13,4.5) -- (14,3.5);
\draw[line width=0.11mm,lightgray] (13,6.5) -- (14,5.5);
\draw[line width=0.11mm,lightgray] (14,5.5) -- (15,4.5);
\draw[line width=0.11mm,lightgray] (12,3.5) -- (13,4.5);
\draw[line width=0.11mm,lightgray] (12,5.5) -- (13,6.5);
\draw[line width=0.11mm,lightgray] (13,4.5) -- (14,5.5);
\draw[line width=0.11mm,lightgray] (14,3.5) -- (15,4.5);
\draw[line width=0.11mm,lightgray] (14,5.5) -- (15,6.5);

\node at (9.7, 3.2) {\tiny{$t_{(0101)}, \, V(y, x+bc)$}};
\draw[line width=0.11mm,lightgray] (11,6.5) -- (11,5.5);
\draw[line width=0.11mm,lightgray] (11,5.5) -- (11,4.5);
\draw[line width=0.11mm,lightgray] (8,3.5) -- (9,3.5);
\draw[line width=0.11mm,lightgray] (9,3.5) -- (10,3.5);
\draw[line width=0.11mm,lightgray] (10,3.5) -- (11,3.5);
\draw[line width=0.11mm,lightgray] (11,3.5) -- (11,4.5);
\draw[line width=0.11mm,lightgray] (9,3.5) -- (9,4.5);
\draw[line width=0.11mm,lightgray] (10,3.5) -- (10,5.5);
\draw[line width=0.11mm,lightgray] (9,4.5) -- (11,4.5);
\draw[line width=0.11mm,lightgray] (9,5.5) -- (9,6.5);
\draw[line width=0.11mm,lightgray] (9,4.5) -- (9,5.5);
\draw[line width=0.11mm,lightgray] (8,5.5) -- (9,5.5);
\draw[line width=0.11mm,lightgray] (9,5.5) -- (10,5.5);
\draw[line width=0.11mm,lightgray] (10,5.5) -- (11,5.5);
\draw[line width=0.2mm] (8,5.5) -- (9,4.5);
\draw[line width=0.2mm] (9,4.5) -- (10,3.5);
\draw[line width=0.2mm] (9,6.5) -- (10,5.5);
\draw[line width=0.2mm] (10,5.5) -- (11,4.5);
\draw[line width=0.11mm,lightgray] (8,3.5) -- (9,4.5);
\draw[line width=0.2mm] (8,5.5) -- (9,6.5);
\draw[line width=0.11mm,lightgray] (9,4.5) -- (10,5.5);
\draw[line width=0.2mm] (10,3.5) -- (11,4.5);
\draw[line width=0.11mm,lightgray] (10,5.5) -- (11,6.5);

\end{tikzpicture}
} 

We collect in the next table the $\mu$-admissible elements and the first three difference vectors for each of them (the fourth is determined by the second and we do not need it). Recall that the finite Weyl group of $\text{GSp}_4$ is isomorphic to $(\Z/2\Z)^2\rtimes S_2$ and it is generated by the elements $s_1$ and $s_2$, corresponding respectively to the generator of $S_2$ and to the second vector of the canonical basis for $(\Z/2\Z)^2$. We are going to use the notation $s_{i_1\cdots i_k}=s_{i_1}\cdots s_{i_k}$.

\begin{center}
\renewcommand{\arraystretch}{1.2}
\begin{tabular}{|p{3cm}|p{3.7cm}|}
\hline
$\text{Adm}(\mu)$ & Difference vectors \\
\hline
$s_{010}\tau=t_{(1100)}$ & $(1100)\;(1100)\;(1100)$ \\ 
$s_{102}\tau=t_{(0101)}$ & $(0101)\;(0101)\;(0101)$ \\
$s_{201}\tau=t_{(1010)}$ & $(1010)\;(1010)\;(1010)$ \\
$s_{212}\tau=t_{(0011)}$ & $(0011)\;(0011)\;(0011)$ \\
$s_{01}\tau=t_{(1100)}s_{2}$ & $(1100)\;(1100)\;(1010)$ \\
$s_{12}\tau=t_{(0101)}s_{2}$ & $(0101)\;(0101)\;(0011)$ \\
$s_{10}\tau=t_{(1100)}s_{121}$ & $(1100)\;(0101)\;(0101)$ \\
$s_{02}\tau=t_{(1010)}s_{1}$ & $(1010)\;(0110)\;(1010)$ \\
$s_{21}\tau=t_{(1010)}s_{121}$ & $(1010)\;(0011)\;(0011)$ \\
$s_{0}\tau=t_{(1100)}s_{21}$ & $(1100)\;(0110)\;(1010)$ \\
$s_{1}\tau=t_{(1100)}s_{1212}$ & $(1100)\;(0101)\;(0011)$ \\
$s_{2}\tau=t_{(1010)}s_{12}$ & $(1010)\;(0110)\;(0011)$ \\
$\tau=t_{(1100)}s_{212}$ & $(1100)\;(0110)\;(0011)$ \\
\hline
\end{tabular}
\end{center} 

\subsubsection*{Length 3}

\begin{itemize}
    \item $t_{(1,1,0,0)}=s_{010}\tau$ has associated difference vector $t^{s_{010}\tau}=\{ (1100)\;(1100)\;(1100) \}$, so for each matrix appearing in (\ref{matr 
    equ g=2}) we need to invert its $2\times 2$ minor given by the first two rows. This means inverting the elements $a^1_{22}$ and $a^2_{11}a^2_{22}-a^2_{12}a^2_{21}$, which correspond to $y$ and $(y+ac)y$ in the ring $B$. Therefore we find that:
    \begin{itemize}
        \item $V_{s_{010}\tau}=\text{Spec} \left( \scalebox{1.2}{$\frac{\Zp[x,y^\pm,a,b,c,(y+ac)^{-1}]}{(xy-p, ax+by+abc)}$} \right) \cong \text{Spec}\left( \Zp[y^\pm,a,c,(y+ac)^{-1}] \right)$
        \item $S_{s_{010}\tau}=\text{Spec}\left( \Fp[y^\pm,a,c,(y+ac)^{-1}] \right)$
    \end{itemize}

    \item $s_{102}\tau$. Invert $x,c,(y+ac)$ and get:
    \begin{itemize}
        \item $V_{s_{102}\tau}= \text{Spec} (\Zp[x^\pm,a,c^\pm,(y+ac)^{-1}])$;
        \item $S_{s_{102}\tau}= \text{Spec} (\Fp[x^\pm,a,c^\pm,(y+ac)^{-1}])$.
    \end{itemize}

    \item $s_{201}\tau$. Invert $b,y,a$ and get:
    \begin{itemize}
        \item $V_{s_{102}\tau}= \text{Spec} (\Zp[a^\pm,b^\pm, y^\pm])$;
        \item $S_{s_{102}\tau}= \text{Spec} (\Fp[a^\pm,b^\pm, y^\pm])$.
    \end{itemize}

    \item $s_{212}\tau$. Invert $x,(x+bc)$ and get:
    \begin{itemize}
        \item $V_{s_{102}\tau}= \text{Spec} (\Zp[x^\pm,b,c, (x+bc)^{-1}])$;
        \item $S_{s_{102}\tau}= \text{Spec} (\Fp[x^\pm,b,c, (x+bc)^{-1}])$.
    \end{itemize}
\end{itemize}

\subsubsection*{Length 2}

\begin{itemize}
    \item $s_{01}\tau$. Invert $y,a$ and get:
    \begin{itemize}
        \item $V_{s_{01}\tau}=\text{Spec}\left( \scalebox{1.2}{$\frac{\Zp[x,y^\pm,a^\pm,b,c]}{(xy-p,ax+by+abc)}$} \right) \cong \text{Spec}\left( \scalebox{1.2}{$\frac{\Zp[y^\pm,a^\pm,b,c]}{(bc-p)}$} \right), \quad
        \left.
        \begin{array}{l}
        x\\[-4pt]
        b \\[-4pt]
        c
        \end{array}
        \right\}
        \mapsto
        \left\{
        \begin{array}{l}
        p/y\\[-4pt]
        -ab/y^2\\[-4pt]
        y(c-1)/a
        \end{array}
        \right.$;

        \item $S_{s_{01}\tau}=(V_{s_{01}\tau})_{\Fp}^\text{sing}=\text{Spec}(\Fp[a^\pm,y^\pm]) \quad (b=c=0)$. 
    \end{itemize}

    \item $s_{12}\tau$. Invert $x,c$ and get:
    \begin{itemize}
        \item $V_{s_{12}\tau}=\text{Spec}\left( \scalebox{1.2}{$\frac{\Zp[x^\pm,y,a,b,c^\pm]}{(xy-p,ax+by+abc)}$} \right) \cong \text{Spec}\left( \scalebox{1.2}{$\frac{\Zp[x^\pm,a,b,c^\pm]}{(a+abc+bp)}$} \right), \quad
        \left.
        \begin{array}{l}
        y\\[-4pt]
        a \\[-4pt]
        b
        \end{array}
        \right\}
        \mapsto
        \left\{
        \begin{array}{l}
        p/x\\[-4pt]
        a/x\\[-4pt]
        bx
        \end{array}
        \right.$;

        \item $S_{s_{12}\tau}=(V_{s_{01}\tau})_{\Fp}^\text{sing}=\text{Spec}(\Fp[x^\pm,c^\pm]), \quad (a=0, b=-1/c)$.
    \end{itemize}
    Note that $S_{s_{12}\tau}$ is contained in the locus where $1+bc=0$, which lies inside the locus where $b$ is invertible. Thus we can replace $V_{s_{12}\tau}$ with
    \begin{gather*}
    V'_{s_{12}\tau}=\text{Spec}\left( \frac{\Zp[x^\pm,a,b^\pm,c^\pm]}{(a+abc+bp)} \right) \cong \text{Spec}\left( \frac{\Zp[x^\pm,a,b^\pm,c, (c-b^{-1})^{-1}]}{(ac-p)} \right)\\
        \quad \;\: \left.
        \begin{array}{l}
        a \\[-4pt]
        c
        \end{array}
        \right\}
        \mapsto
        \left\{
        \begin{array}{l}
        -a\\[-4pt]
        c-b^{-1}
        \end{array}
        \right. .
    \end{gather*}

    \item $s_{10}\tau$. Invert $c, y+ac$ and get:
    \begin{itemize}
        \item 
        $V_{s_{10}\tau}=\text{Spec}\left( \scalebox{1.2}{$\frac{\Zp[x,y,a,b,c^\pm, (y+ac)^{-1}]}{(xy-p,ax+by+abc)}$} \right) \cong \text{Spec}\left( \scalebox{1.2}{$\frac{\Zp[x,y,a,c^\pm, (y+ac)^{-1}]}{(xy-p)}$} \right), \\  b\mapsto -\frac{ax}{(y+ac)};$

        \item $S_{s_{10}\tau}=(V_{s_{01}\tau})_{\Fp}^\text{sing}=\text{Spec}(\Fp[a^\pm,c^\pm]), \quad (y=x=0).$
    \end{itemize}

    \item $s_{02}\tau$. Invert $a,b$ and get:
    \begin{itemize}
        \item $V_{s_{01}\tau}=\text{Spec}\left( \scalebox{1.2}{$\frac{\Zp[x,y,a^\pm,b^\pm,c]}{(xy-p,ax+by+abc)}$} \right) \cong \text{Spec}\left( \scalebox{1.2}{$\frac{\Zp[a^\pm,b^\pm,x,c]}{(xc-p)}$} \right), \quad
        \left.
        \begin{array}{l}
        y \\
        c
        \end{array}
        \right\}
        \mapsto
        \left\{
        \begin{array}{l}
        c\\
        -\frac{cb+ax}{ab}
        \end{array}
        \right.$;

        \item $S_{s_{01}\tau}=(V_{s_{01}\tau})_{\Fp}^\text{sing}=\text{Spec}(\Fp[a^\pm,b^\pm]), \quad (x=c=0).$
    \end{itemize}

    \item $s_{21}\tau$. Invert $b, x+bc$ and get:
    \begin{itemize}
        \item $V_{s_{01}\tau}=\text{Spec}\left( \scalebox{1.2}{$\frac{\Zp[x,y^\pm,a^\pm,b,c]}{(xy-p,ax+by+abc)}$} \right) \cong \text{Spec}\left( \scalebox{1.2}{$\frac{\Zp[y^\pm,a^\pm,b,c]}{(xy-p)}$} \right), \quad a\mapsto -by/(x+bc)$;

        \item $S_{s_{01}\tau}=(V_{s_{01}\tau})_{\Fp}^\text{sing}=\text{Spec}(\Fp[a^\pm,y^\pm]), \quad (x=y=0).$
    \end{itemize}
\end{itemize}

\subsubsection*{Length 1}

\begin{itemize}
    \item $s_1\tau$. The stratum $S_{s_1\tau}$ lies in the intersection of all four irreducible components of the special fiber, which is given by 
    \[
    V(x,y,a,b) \subset \text{Spec}\left( \frac{\Fp[x,y,a,b,c]}{(xy, ax+by+abc)} \right).
    \]
    Removing the worst singular point from this intersection we obtain the equations of the stratum, which is therefore
    \[
    S_{s_1\tau}=\text{Spec}(\Fp[c^\pm]) \xrightarrow{\makebox[1cm]{\tiny{$\substack{x=y=0\\a=b=0}$}}} \text{Spec}\left( \frac{\Fp[x,y,a,b,c]}{(xy, ax+by+abc)} \right).
    \]
    
    \item $s_0\tau$. The stratum $S_{s_0\tau}$ lies in the intersection of the three irreducible components $V(x,b), V(x, y+ac), V(y, x+bc)$ which is given by
    \[
    \text{Spec}\left( \frac{\Fp[a,c]}{(ac)}\right) \xrightarrow{\makebox[1cm]{\tiny{$\substack{x=y=b=0\\ac=0}$}}} \text{Spec}\left( \frac{\Fp[x,y,a,b,c]}{(xy, ax+by+abc)} \right).
    \]
    From this, we need to remove $S_{s_1\tau}\cup \tau$, which corresponds to $a=0$. Thus we get that
    \[
    S_{s_0\tau}=\text{Spec}(\Fp[a^\pm]) \xrightarrow{\makebox[1cm]{\tiny{$\substack{x=y=0\\c=b=0}$}}} \text{Spec}\left( \frac{\Fp[x,y,a,b,c]}{(xy, ax+by+abc)} \right).
    \]
    
    \item $s_2\tau$. The stratum $S_{s_2\tau}$ lies in the intersection of the irreducible components $V(x, y+ac), V(y,a), V(y, x+bc)$. Proceeding as in the previous case we find that
    \[
    S_{s_2\tau}=\text{Spec}(\Fp[b^\pm])\xrightarrow{\makebox[1cm]{\tiny{$\substack{x=y=0\\a=c=0}$}}} \text{Spec}\left( \frac{\Fp[x,y,a,b,c]}{(xy, ax+by+abc)} \right).
    \]
\end{itemize}

\subsubsection*{Length 0}
There is a unique length 0 admissible element $\tau$, whose open containing it is the whole $U$ and whose stratum $S_\tau$ consists of the worst singular point given by $x=y=a=b=c=0$.

\subsection{Pro-\texorpdfstring{$p$}{} Iwahori case}\label{resolution pro p}

We can now construct a semistable resolution of the local model for the Siegel modular variety with pro-$p$ Iwahori level structure in the case of $\text{GSp}_4$. We explain how to do this following Shadrach's approach as in \cite[Chapter 6]{Sha}, where he performs a series of blow-ups with centers in the special fiber. We also describe a modification to his procedure which is shorter and hence better suited for trace of Frobenius computations. The $\Gamma_1(p)$-integral model is the spectrum of the ring
\[
C=\frac{B[u_0,u_1,v_0,v_1]}{\left( \scalebox{1.2}{$\substack{u_0^{p-1}-y,\, v_0^{p-1}-x,\, u_1^{p-1}-(y+ac)\\ v_1^{p-1}-(x+bc),\, u_0v_0-u_1v_1}$} \right)}.
\]

\subsubsection*{Step 1: base change to $U'$.}

Consider the fiber product
\[
U_1'=U_1\times_{U}U'=\text{Proj}_B\left( \frac{B[u_i,v_i][\tilde{x}, \tilde{b}]}{\left( \scalebox{1.2}{$\substack{u_0^{p-1}-y,\; v_0^{p-1}-x,\; u_1^{p-1}-(y+ac), \; v_1^{p-1}-(x+bc)\\ u_0v_0-u_1v_1,\; x\tilde{b}-b\tilde{x},\; a\tilde{x}+\tilde{b}y+a\tilde{b}c}$} \right)} \right),
\]
with $\tilde{x}, \tilde{b}$ of degree 1 and the other variables of degree 0; call $\sigma_1:U_1'\longrightarrow U_1$ the projection. Restricting to the two charts of $U'$ we find the following equations:
\begin{itemize}
    \item[] 
    \[
    U'_{1,[\tilde{b}=1]}=U_1'\times_{U_0'}U'_{[\tilde{b}=1]}=\text{Spec}\left(\frac{\Zp[\tilde{x},a,b,c]}{(\tilde{x}abc-p)}[u_i,v_i]\Big/\left( \scalebox{1.2}{$\substack{u_0^{p-1}-ac,\; u_1^{p-1}-a\tilde{x},\; v_0^{p-1}-b\tilde{x}\\ v_1^{p-1}-bc,\; u_0v_0-u_1v_1}$} \right)\right);
    \]

    \item[]
    \[
    U'_{1,[\tilde{x}=1]}=U_1'\times_{U_0'}U'_{[\tilde{x}=1]}=\text{Spec}\left(\frac{\Zp[x,y,a,\tilde{b},c]}{\left(\scalebox{1.2}{$\substack{a+\tilde{b}y+a\tilde{b}c\\xy-p}$}\right)}[u_i,v_i] \Big/ \left( \scalebox{1.2}{$\substack{u_0^{p-1}-y,\; v_0^{p-1}-x,\; u_1^{p-1}-(y+ac)\\ \; v_1^{p-1}-x(1+\tilde{b}c),\; u_0v_0-u_1v_1}$} \right)\right) .
    \]
\end{itemize}

\begin{remark}\label{different resolution}
Shadrach proceeds with further blow-ups, but we can find directly a semistable resolution of the two opens ${D}(\tilde{b})$ and ${D}(1+\tilde{b}c)$ covering the $[\tilde{x}=1]$-chart: this will come very handy in the computations with nearby cycles.
\begin{itemize}
    \item On $D(1+\tilde{b}c)$ we get 
    \[
    U'_{1,[\tilde{x}=1, (1+\tilde{b}c)^{-1}]}=\text{Spec}\left(\frac{\Zp[x,y,\tilde{b},c,(1+\tilde{b}c)^{-1}]}{(xy-p)}[u_i,v_i]\Big/\left( \scalebox{1.2}{$\substack{u_0^{p-1}-y,\; u_1^{p-1}-y/(1+\tilde{b}c),\; v_0^{p-1}-x\\ v_1^{p-1}-x(1+\tilde{b}c),\; u_0v_0-u_1v_1}$} \right)\right),
    \]
    which admits a map from
    \[
    \text{Spec}\left(\frac{\Zp[x,y,\tilde{b},c,(1+\tilde{b}c)^{-1}]}{(xy-p)}[r,s,t]\Big/\left( \scalebox{1.2}{$\substack{r^{p-1}-x,\;, s^{p-1}-y\\ t^{p-1}-(1+\tilde{b}c)}$} \right)\right)
    \]
    by sending $u_0,u_1,v_0,v_1$ respectively to $r,rt,s,s/t$. This is a finite map, which has an inverse over $\Qp$ given by $r, s, t \longmapsto u_0, v_0, u_1/u_0$.

    \item On $D(\tilde{b})$ instead we get
    \[
    U'_{1,\left[{\tilde{x}=1, \tilde{b}^{-1}}\right]}=\text{Spec}\left(\frac{\Zp[x,a,\tilde{b}^\pm,c]}{(xac-p)}[u_i,v_i]\Big/\left( \scalebox{1.2}{$\substack{u_0^{p-1}-ac,\; u_1^{p-1}-a\tilde{b}^{-1},\; v_0^{p-1}-x\\ v_1^{p-1}-\tilde{b}xc,\; u_0v_0-u_1v_1}$} \right)\right),
    \]
    which can be mapped from
    \[
    \text{Spec}\left(\frac{\Zp[x,a,\tilde{b}^\pm,c]}{(xac-p)}[r,s,t]\Big/\left( \scalebox{1.2}{$\substack{r^{p-1}-x,\; s^{p-1}-a\tilde{b}^{-1}\\ t^{p-1}-c\tilde{b}}$} \right)\right),
    \]
    by sending $u_0,u_1,v_0,v_1$ to $r,rt,st,s$. Again, this is a finite map with an inverse over $\Qp$.
\end{itemize}
\end{remark}

\subsubsection*{Step 2: blow-up along $(v_0,v_1)$.}

We now consider $U_1''=\text{Bl}_{(v_0,v_1)}(U_1')$, with projection $\sigma_2:U_1''\longrightarrow U_1'$. It is not hard to check, by working with its two charts, that it can be described as 
\[
\text{Proj}_{U_1'}\left( \mathcal{O}_{U_1'}[\tilde{v}_0,\tilde{v}_1] \Big/ \left( \scalebox{1.2}{$\substack{v_0\tilde{v}_1-v_1\tilde{v}_0,\; u_0\tilde{v}_0-u_1\tilde{v}_1\\ \tilde{v}_0^{p-1}(\tilde{x}+\tilde{b}c)-\tilde{v}_1^{p-1}\tilde{x}}$} \right) \right).
\]
Restricting to the two opens $U'_{1,[\tilde{b}=1]}, U'_{1,[\tilde{x}=1]}$ of $U_1'$ we find:
\begin{itemize}
    \item[]
    \[
    U''_{1,[\tilde{b}=1]}=\text{Proj}_{U'_{1,[\tilde{b}=1]}}\left( \mathcal{O}_{U'_{1,[\tilde{b}=1]}}[\tilde{v}_0,\tilde{v}_1]\Big/\left( \scalebox{1.2}{$\substack{v_0\tilde{v}_1-v_1\tilde{v}_0,\; u_0\tilde{v}_0-u_1\tilde{v}_1\\ \tilde{x}\tilde{v}_1^{p-1}-(\tilde{x}+c)\tilde{v}_0^{p-1}}$} \right) \right),
    \]

    \item[]
    \[
    U''_{1,[\tilde{x}=1]}=\text{Proj}_{U'_{1,[\tilde{x}=1]}}\left( \mathcal{O}_{U'_{1,[\tilde{x}=1]}}[\tilde{v}_0,\tilde{v}_1]\Big/\left( \scalebox{1.2}{$\substack{v_0\tilde{v}_1-v_1\tilde{v}_0,\; u_0\tilde{v}_0-u_1\tilde{v}_1\\ \tilde{v}_1^{p-1}-(1+\tilde{b}c)\tilde{v}_0^{p-1}}$} \right) \right).
    \]
\end{itemize}

After some calculations, one can check that the $[\tilde{v}_0=1]$-chart of $U''_{1,[\tilde{b}=1]}$, respectively the $[\tilde{v}_1=1]$-chart, is isomorphic to
    \[
    U''_{1,[\tilde{b}=\tilde{v}_0=1]}=\text{Spec}\left(\frac{\Zp[u_1,v_0,\tilde{v}_1,a,b,\tilde{x}]}{\left( \scalebox{1.2}{$\substack{u_1^{p-1}-a\tilde{x},\; v_0^{p-1}-b\tilde{x}\\ \tilde{v}_1^{p-1}ab\tilde{x}^2-p}$} \right)}\right), \quad U''_{1,[\tilde{b}=\tilde{v}_1=1]}=\text{Spec}\left(\frac{\Zp[u_0,v_1,\tilde{v}_0,a,b,c]}{\left( \scalebox{1.2}{$\substack{u_0^{p-1}-ac,\; v_1^{p-1}-bc\\ \tilde{v}_0^{p-1}abc^2-p}$} \right)}\right).
    \]
Similarly, the $[\tilde{v}_0=1]$-chart of $U''_{1,[\tilde{x}=1]}$, respectively the $[\tilde{v}_1=1]$-chart, is isomorphic to
    \[
    U''_{1,[\tilde{x}=\tilde{v}_0=1]}=\text{Spec}\left(\frac{\Zp[\tilde{b},c,u_1,v_0,\tilde{v}_1]}{\left( \scalebox{1.2}{$\substack{(u_1v_0\tilde{v}_1)^{p-1}-p\\ \tilde{v}^{p-1}_1-(1+\tilde{b}c)}$} \right)}\right), \quad   U''_{1,[\tilde{x}=\tilde{v}_1=1]}=\text{Spec}\left(\frac{\Zp[\tilde{b},c,u_0,v_1,\tilde{v}_0^\pm]}{\left( \scalebox{1.2}{$\substack{(u_0v_1\tilde{v}_0)^{p-1}-p\\ \tilde{v}^{p-1}_0-(1+\tilde{b}c)}$} \right)}\right).
    \]
Note that $U''_{1,[\tilde{x}=\tilde{v}_1=1]}$ is an étale cover of a scheme with semistable reduction, so it has semistable reduction. The same holds for $U''_{1,[\tilde{x}=\tilde{v}_0=1]}$ restricted to the locus where $1+\tilde{b}c$ is invertible, whereas on the locus where $\tilde{b}$ is invertible we can apply the change of variables $c\mapsto \tilde{b}^{-1}(c-1)$ and directly see that it has semistable reduction.

\subsubsection*{Step 3: blow-up along $(u_0,u_1,v_0,v_1)$.}

We now blow-up $U_1''$ along the closed subscheme cut out by the ideal $(u_0,u_1,v_0,v_1)$, which on each of the four charts from step 2 becomes:
\begin{align*}
    &(u_1,v_0) \;\; \text{on}\;\; U''_{1,[\tilde{b}=\tilde{v}_0=1]}, &(u_0, v_1) \;\; \text{on}\;\; U''_{1,[\tilde{b}=\tilde{v}_1=1]},\\
    &(u_1, v_0) \;\; \text{on}\;\; U''_{1,[\tilde{x}=\tilde{v}_0=1]}, &(u_0,v_1) \;\; \text{on}\;\; U''_{1,[\tilde{x}=\tilde{v}_1=1]}.
\end{align*}
The resulting scheme $U_1'''=\text{Bl}_{(u_i,v_i)_i}(U_1'')\xrightarrow[]{\sigma_3}U_1''$ is covered by eight charts. The four charts coming from $U'_{1,[\tilde{b}=1]}$ are all isomorphic to the spectrum of
\[
E=\frac{\Zp[r,s,t,e,f]}{\left( \scalebox{1.2}{$\substack{(r^2st)^{p-1}-p\\ r^{p-1}-ef}$}\right) }.
\]
The four charts coming from $U'_{1,[\tilde{x}=1]}$ remain almost untouched, except that in each case one of the variables $u_i$ or $v_i$ appears in the equations with a power of 2.

\subsubsection*{Step 4: iterated blow-ups along $(u_i,v_i,\tilde{x},c)_i$.}

We now blow-up $U_1'''$ along the ideal $(u_0,u_1,v_0,v_1,\tilde{x},c)$, which does not intersect the $[\tilde{x}=1]$-charts, so they will not be modified. On each of the four charts isomorphic to $\text{Spec}(E)$ this ideal corresponds to the ideal $(r,e)$ and 
\[
\text{Bl}_{(r,e)}(\text{Spec}(E))=\text{Proj}\left( E[\tilde{r},\tilde{e}] \Big/ \left( \scalebox{1.2}{$\substack{\tilde{r}e-\tilde{e}r \\ \tilde{r}r^{p-2}-\tilde{e}f}$} \right) \right)
\]
with the two charts given by
\[
\text{Spec}\left( \frac{\Zp[r,s,t,\tilde{e},f]}{\left( \scalebox{1.2}{$\substack{(r^2st)^{p-1}-p\\ r^{p-2}-\tilde{e}f}$} \right)} \right) \quad \text{and} \quad \text{Spec}\left( \frac{\Zp[\tilde{r},s,t,e,f]}{\left( {(\tilde{r}^2e^2st)^{p-1}-p} \right)} \right).
\]
The second chart is already semistable, while in the first the exponent of $r$ has dropped by 1: we can proceed inductively to obtain a semistable resolution of $U_1$. More precisely, for $i=0,\ldots,p-2$ let
\[
E_i=\frac{\Zp[r,s,t,e_i,f]}{\left( \scalebox{1.2}{$\substack{(r^2st)^{p-1}-p\\ r^{p-1-i}-e_if}$} \right)} \quad \text{and} \quad R_{i-1}=\frac{\Zp[\tilde{r},s,t,e_{i-1},f]}{\left( {(\tilde{r}^2e_{i-1}^2st)^{p-1}-p} \right)},
\]
so that for $i=1,\ldots,p-2$ we have
\[
\text{Bl}_{(r,e_{i-1})}(\text{Spec}(E_{i-1}))=\text{Proj}\left( E_{i-1}[\tilde{r},e_i] \Big/ \left( \scalebox{1.2}{$\substack{\tilde{r}e_{i-1}-e_ir \\ \tilde{r}r^{p-1-i}-e_if}$} \right) \right)=\text{Spec}(E_i)\cup \text{Spec}(R_{i-1}):
\]
at the $(p-2)$-th step we obtain that $E_{p-2}\cong R_{p-3}$, which is semistable, so the procedure is complete.

\section{Trace of Frobenius on nearby cycles on the local model}\label{nearby cycles}

In this section we consider the action of the arithmetic Frobenius at $q$ on the complex of nearby cycles, which we briefly introduce in the first subsection. We are going to compute its trace on the local model $U_1$ using the resolution we just constructed.

\subsection{Nearby cycles}

Fix a henselian triple $(Z,s,\eta)$ with $s$ of characteristic $p$, i.e. $Z$ is the spectrum of a henselian discrete valuation ring, with special and generic point $s$ and $\eta$ respectively: the relevant example for us will be $(\Zp, \Fp, \Qp)$. Choose an algebraic closure $\overline{\eta}$ of $\eta$, which determines the normalization $\overline{Z}$ of $Z$ with residue field $\overline{s}$. For a scheme $X$ over $Z$, denote by $X_s$, $X_\eta$ its special and generic fiber, by $X_{\overline{Z}}$ its base change to $\overline{Z}$, with special and generic fiber $X_{\overline{s}}$, $X_{\overline{\eta}}$. They fit in the following diagrams, where all squares are cartesian:
\[
\begin{tikzcd}
    X_s \arrow[d] \arrow[r, "i", hook] & X \arrow[d] & X_\eta \arrow[d] \arrow[l, "j"', hook'] & X_{\overline{s}} \arrow[d] \arrow[r, "\bar{i}", hook] & X_{\overline{Z}} \arrow[d] & X_{\overline{\eta}} \arrow[d] \arrow[l, "\bar{j}"', hook'] \\
    s \arrow[r, hook] & Z & \eta \arrow[l, hook'] & \overline{s} \arrow[r, hook] & \overline{Z} & \overline{\eta} \arrow[l, hook']     
\end{tikzcd}
\]
Let $D^b_c(X, \overline{\Q}_\ell)$ be the derived category of $\overline{\Q}_\ell$-sheaves on $X$, $D^b_c(X \times_s \eta, \overline{\Q}_\ell)$ the category of objects $\mathcal{F}\in D^b_c(X_{\overline{s}}, \overline{\Q}_\ell)$ together with a continuous $\text{Gal}(\overline{\eta}/\eta)$-action, compatible with the action on $X_{\bar{s}}$ (continuity is tested on cohomology sheaves). For $\mathcal{F}\in D^b_c(X_{\eta}, \overline{\Q}_\ell)$, the complex of nearby cycles of $\mathcal{F}$ is
\[
R\Psi(\mathcal{F})=\bar{i}^*R\bar{j}_*(\mathcal{F}_{\overline{\eta}}),
\]
where $\mathcal{F}_{\overline{\eta}}$ is the pullback of $\mathcal{F}$ to $X_{\overline{\eta}}$: it is an object of $D^b_c(X \times_s \eta, \overline{\Q}_\ell)$. When $X$ has semistable reduction with all multiplicities of the divisors of the special fiber prime to $p$, the nearby cycles complex for the constant sheaf $\overline{\Q}_\ell$ can be described quite explicitly: we refer the reader to \cite[Chapter 5]{HR} for the following discussion and to \cite[Section 2.1]{Del} for the original source on this material. Let $R^i\Psi=\bar{i}^*R^i\bar{j}_*(\overline{\Q}_\ell)$, for a point $x\in X_{\overline{s}}$ choose a geometric point $\bar{x}$ over it and let $\mathcal{C}_x=\{C_i\}_i$ be the set of branches of $X_{\overline{s}}$ passing through it, with multiplicity $d_i$. Associated to $x$, we have the following homological complex concentrated in degree 1 and 0:
\begin{align*}
    K_x:&\Z^\mathcal{C_x}\xrightarrow{\makebox[1cm]{}} \Z\\
    (&n_i)_i\longmapsto \sum_id_in_i.
\end{align*}
Since the $d_i$ are all prime to $p$, the action of $\text{Gal}(\overline{\eta}/\eta)$ factors through the tame fundamental group $\text{Gal}(\eta_t/\eta)$ and we have Galois-equivariant isomorphisms
\begin{align*}
    &R^0\Psi_{\bar{x}}=\overline{\Q}_\ell^{H_0(K_x)},\\
    &R^i\Psi_{\bar{x}}=\bigwedge^i(H_1(K_x)\otimes \overline{\Q}_\ell(-1))\otimes R^0\Psi_{\bar{x}}.
\end{align*}
Here, the action of the tame inertia $I_t\subset \text{Gal}(\eta_t/\eta)$ on $\bigwedge^i(H_1(K_x)\otimes \overline{\Q}_\ell(-1))$ is trivial and the Frobenius acts via its natural action on $\overline{\Q}_\ell(-1)$. On the other hand, $I_t$ acts on $R^0\Psi_{\bar{x}}$ by permuting transitively the elements of $H_0(K_x)$. In fact, this set can be seen as the set of connected components of the strict henselization of $X$ in $\bar{x}$: it is easy to check that its connected components are indexed by the $g$-th roots of unity in $\overline{\Z}_p$, where $g=\text{gcd}(d_i)_i$, and that the tame inertia acts on them through its natural quotient $\mu_g$ by multiplication. In particular, the inertia invariants on the stalk of the $i$-th cohomology sheaf of the nearby cycles are
\[
(R^i\Psi_{\bar{x}})^I=\bigwedge^i(H_1(K_x)\otimes \overline{\Q}_\ell(-1)),
\]
which implies that the alternating trace of the Frobenius on inertia invariants of $R\Psi_{\bar{x}}$ is given by
\[
\text{Tr}(\text{Fr}_q; \,(R\Psi_{\bar{x}})^I)=\sum_i(-1)^i\text{Tr}(\text{Fr}_q; \,(R^i\Psi_{\bar{x}})^I)=(1-q)^{\#\mathcal{C}_x-1}
\]
where $q$ is the cardinality of the residue field of $x$. Of course, the local model $U_1$ is not semistable, so we cannot directly apply this. However, we have the following proposition, which is a consequence of the base change theorems for proper or smooth morphisms.
\begin{proposition}[\cite{Del} Section 2.1.7]
Let $f:X'\longrightarrow X$ be a map of schemes over $Z$ and denote by $R\Psi$, $R\Psi'$ the complex of nearby cycles on $X$, $X'$ respectively.
    \begin{itemize}
        \item[$(i)$] If $f$ is proper, the natural morphism
            \[
            R\Psi Rf_{\eta,*}\longrightarrow Rf_{s,*}R\Psi'
            \]
            is an isomorphism.
        \item[$(ii)$] If $f$ is smooth, the natural morphism
            \[
            f_s^*R\Psi \longrightarrow R\Psi'f_\eta^*
            \]
            is an isomorphism.
    \end{itemize}
\end{proposition}
In particular, if $f$ is proper and $f_\eta$ is an isomorphism we find that 
\begin{equation}\label{stalk of pushforward}
R\Psi_{\bar{x}}=\bigoplus_{\bar{x}'\overset{f_s}{\longmapsto}\bar{x}} R\Psi'_{\bar{x}'},
\end{equation}
where $x\in X(\mathbb{F}_q)$ and the sum is indexed by geometric points of $X'$ lying above $\bar{x}$: this decomposition is preserved after passing to inertia invariants. If all such points lie in $X'(\mathbb{F}_q)$ we obtain that
\[
\text{Tr}(\text{Fr}_q; (R\Psi_{\bar{x}})^I)=\sum_{\bar{x}'\overset{f_s}{\longmapsto}\bar{x}}\text{Tr}(\text{Fr}_q; (R\Psi'_{\bar{x}'})^I).
\]
However, it can happen that some of these points are only defined over a proper extension $\mathbb{F}_{q'}$ of $\mathbb{F}_q$: in this case the action of the Frobenius map on the direct sum (\ref{stalk of pushforward}) is no longer diagonal, so the above formula for the trace is not valid anymore. We are going to run into the situation where the whole fiber $f_s^{-1}(\bar{x})$ lies outside of $X'(\mathbb{F}_q)$ and $\text{Fr}_q$ permutes it transitively: in this case we can conclude that the trace of Frobenius on $(R\Psi_{\bar{x}})^I$ is 0. In the non-semistable situation we do not know a priori that the inertia action is exact, so we are led to work with the semi-simple trace of Frobenius on nearby cycles, for which we refer to \cite{HN} or \cite[Section 9.3]{Hai}. From now on, we will write $\tr(\fr_q,R\Psi_x)$ for the semi-simple trace of Frobenius on $R\Psi_x$. We will also avoid to say semi-simple every time, but when we talk about the trace of Frobenius it should be intended as the semi-simple trace.

\begin{notation}
    If we need to specify that we are considering nearby cycles on a scheme $X$, we will write $R\Psi^X$. For $i=0,1$, we are going to write $R\Psi^i$ for the nearby cycles on $\mathcal{A}_i$.
\end{notation}

\subsection{Trace of Frobenius on \texorpdfstring{$U_1$}{}}\label{trace of frob on U1}
We can now compute the trace of Frobenius on nearby cycles on the local model $U_1$: for ease of notation, we will write $R\Psi$ instead of $R\Psi^{U_1}$. We are going to perform our computations on the preimage of each stratum $S_w\subset U_{\Fp}$ along the map
\[
\pi:U_1\longrightarrow U
\]
separately. Moreover, we can restrict to $\pi^{-1}(V_w)$: this is an advantage, since $V_w$ is often semistable and it is easier to resolve its preimage. Let us write $V_{1,w}=\pi^{-1}(V_w)$ and $S_{1,w}=\pi^{-1}(S_w)$ for $w\in\text{Adm}(\mu)$. It turns out that in all except two cases the trace will be constant on $S_{1,w}$ and equal to the trace on $S_w$.

\subsubsection*{Length 3}

\begin{itemize}
    \item $s_{010}\tau$. The open $V_{1,s_{010}\tau}$ is isomorphic to the spectrum of
    \[
     \frac{\Zp[y^\pm,a,c,(y+ac)^{-1}][u_0,v_0,v_1]}{\left( \scalebox{1.2}{$\substack{v_0^{p-1}-y, \, v_1^{p-1}-(y+ac) \\ u_0^{p-1}-p }$} \right)}, 
    \]
    thus for any $x\in S_{1,s_{010}\tau}(\mathbb{F}_q)$ we get that $\text{Tr}(\text{Fr}_q; R\Psi_{{x}})=1$. 

    \item The other three open strata work exactly the same: the preimage of the corresponding open $V$ is of the same shape as above (two $(p-1)$-th roots of units and one $(p-1)$-th root of $p$) and the trace of Frobenius is constantly 1 on the preimage of the stratum.
    
\end{itemize}

\subsubsection*{Length 2}

\begin{itemize}
    \item $s_{01}\tau$. The open $V_{1,s_{01}\tau}$ is isomorphic to the spectrum of
    \[
    \frac{\Zp[y^\pm, a^\pm, u_0,u_1,v_1]}{\left( (u_1v_1)^{p-1}-p, u_0^{p-1}-y \right)}
    \]
    and $S_{1,s_{01}\tau}$ is the locus where $u_1^{p-1}=v_1^{p-1}=0$. Therefore for any point $x\in S_{1,s_{01}\tau}(\mathbb{F}_q)$ we get that $\text{Tr}(\text{Fr}_q; R\Psi_{{x}})=1-q$.
    
    \item As in the length 3 case, for the strata corresponding to $s_{10}\tau, s_{21}\tau, s_{12}\tau$ the calculations are completely analogous and one gets the same result (using the modified open $V'_{s_{12}\tau}$ for $s_{12}\tau$).

    \item $s_{02}\tau$. In this case the situation is different, so we describe it in some more details. This is not surprising, since the stratum corresponding to $s_{02}\tau$ is containted in the $p$-rank zero locus of $\AI$, see \cite[Example 2.6]{GY}. First, the open $V_{1,s_{02}\tau}$ is isomorphic to the spectrum of the ring
    \[
    A=\frac{\Zp[a^\pm, b^\pm, x, c]}{\left( xc-p \right)}[u_0,u_1,v_0,v_1] \Big/ \left( \scalebox{1.2}{$\substack{ u_0^{p-1}-c,\, u_1^{p-1}-ax/b\\ v_0^{p-1}-x,\, v_1^{p-1}-cb/a\\ u_0v_0-u_1v_1}$} \right),
    \]
    with $S_{1,s_{02}\tau}$ corresponding to $x=c=0$. Thus the $\mathbb{F}_q$-points of $S_{1,s_{02}\tau}$ are given by setting to 0 all the variables of $A$ except for $a$ and $b$, which can take any value in $\Fq^\times$: we denote by $(\alpha, \beta)$ the $\mathbb{F}_q$-point of $S_{1,s_{02}\tau}$ determined by the choice of $\alpha$ and $\beta$ in $\Fq^\times$. We can resolve $\text{Spec}(A)$ with the map
    \[
    A\longrightarrow B=\frac{\Zp[a^\pm, b^\pm, x, c, u_0, v_0, t]}{\left( \scalebox{1.2}{$\substack{xc-p,\, u_0^{p-1}-c\\ v_0^{p-1}-x,\, t^{p-1}-a/b}$} \right)}, \quad u_0,v_0,u_1,v_1\longmapsto u_0,v_0,v_0t,u_0t^{-1},
    \]
    which is finite and admits an inverse over $\Qp$ ($t\longmapsto u_1u_0^{-1}$): we denote by $f$ the induced map at the level of spectra. Its fiber over $(\alpha,\beta)$ consists of points where $a=\alpha, b=\beta, x=c=u_0=v_0=0$ and $t^{p-1}=\alpha/\beta$. There are always $p-1$ such points and they are defined over $\Fq$ if and only if $\alpha/\beta$ is a $(p-1)$-th power in $\Fq^\times$, or equivalently if $(\alpha/\beta)^{\frac{q-1}{p-1}}=1$: if not, they are defined over $\mathbb{F}_{q^{p-1}}$. Moreover, when the fiber is defined over $\Fq$, the trace of Frobenius on $(\alpha,\beta)$ is equal to $1-q$. Therefore, in light of the discussion at the end of the previous subsection, we see that
    \begin{equation}\label{trace s02}
    \text{Tr}(\text{Fr}_q, R\Psi_{(\alpha, \beta)})=\begin{cases}
        0, &\text{if} \;(\alpha/\beta)^{\frac{q-1}{p-1}}\neq1\\
        (p-1)(1-q), &\text{if}\; (\alpha/\beta)^{\frac{q-1}{p-1}}=1
    \end{cases}
    \end{equation}
\end{itemize}

\subsubsection*{Length 1}

\begin{itemize}
    \item $s_{0}\tau, s_{2}\tau$. For these two strata the calculations are the same, so we explain what happens for $s_{2}\tau$. In this case, $V_{1, s_{2}\tau}$ is the spectrum of
    \[
    A=\frac{\Zp[x,a,b^\pm,c]}{(axc-p)}[u_0,u_1,v_0,v_1] \Big/ \left( \scalebox{1.2}{$\substack{u_0^{p-1}-ac,\, u_1^{p-1}-ax\\ v_0^{p-1}-x, \, v_1^{p-1}-c\\ u_0v_0-u_1v_1}$} \right)
    \]
    and $S_{1, s_{2}\tau}$ corresponds to $x=a=c=0$. Again, we can resolve $\text{Spec}(A)$ with the following finite map which is an isomorphism over $\Qp$:
    \[
    A\longmapsto \frac{\Zp[x,a,b^\pm,c,r,s,t]}{\left( \scalebox{1.2}{$\substack{axc-p,\, r^{p-1}-a\\s^{p-1}-x,\, t^{p-1}-c}$} \right)}, \quad u_0,u_1,v_0,v_1\longmapsto rt,rs,s,t.
    \]
    This time, however, the fiber over any $\Fq$ point $x$ of $S_{1, s_{2}\tau}$ consists of a single point, with trace of Frobenius equal to $(1-q)^2$. Thus $\text{Tr}(\text{Fr}_q; R\Psi_{{x}})=(1-q)^2$ in this case.

    \item $s_{1}\tau$. The open $V_{1,s_1\tau}$ is the spectrum of
    \[
    \frac{\Zp[x,y,a,b,c^\pm]}{({xy-p,\, ax+by+abc})}[u_0,u_1,v_0,v_1]\Big/\left( \scalebox{1.2}{$\substack{u_0^{p-1}-y,\, u_1^{p-1}-(y+ac)\\ v_0^{p-1}-x, \, v_1^{p-1}-(x+bc)\\ u_0v_0-u_1v_1}$} \right),
    \]
    with $S_{,1s_1\tau}$ given by $x=y=a=b=0$, so that its $\Fq$-points are given by setting $u_i=v_i=0$ and $c=\gamma$ for $\gamma\in\Fq^\times$: we will denote by $\gamma$ such a point. We are going to use the first two steps of the semistable resolution constructed in Subsection \ref{resolution pro p}. We will adopt the same notation as used there, with the understanding that we are base changing everything over the locus where $c$ is invertible. Consider the preimage of $\gamma$ along the projection $\sigma_1: U'_{1}\rightarrow U_{1}$: it is given by
    \[
    \sigma_1^{-1}(\gamma)=\text{Proj}(\Fq[\tilde{x}, \tilde{b}])=\begin{cases}
        (1:\lambda)\in U'_{1,[{\tilde{x}=1, \tilde{b}^{-1}}]}, \left( \substack{x=a=u_i=v_i=0\\\tilde{b}=-\lambda,\, c=-\gamma-\lambda^{-1}} \right), \, \lambda\in\Fq^\times \\[1ex]
        (1:0) \in U'_{1,[{\tilde{x}=1, (1+\tilde{b}c)^{-1}}]} \left( \substack{x=a=u_i=v_i=0\\\tilde{b}=0,\, c=\gamma} \right)\\[1ex]
        (0:1) \in U'_{1,[\tilde{b}=1]} \left( \substack{a=b=\tilde{x}=0,\, c=\gamma \\ u_i=v_i=0 } \right).
    \end{cases}
    \]
    We compute case by case the trace of Frobenius on each point of this fiber and then sum them up.
    \begin{itemize}
        \item[$(i)$] In the first case, we can use the resolution of Remark \ref{different resolution}: one checks that the fiber over $(1:\lambda)$ is given by $c=-\gamma-\lambda^{-1}, \tilde{b}=-\lambda, t^{p-1}=1+\gamma\lambda$ and all the other variables equal to 0. There are three possible cases:
        \begin{itemize}
            \item $\lambda=-\gamma^{-1}$: only one point over $(1:\lambda)$, trace equals to $(1-q)^2$;
            \item $(1+\lambda\gamma)^{\frac{q-1}{p-1}}=1$: $(p-1)$ $\Fq$-points over $(1:\lambda)$, trace equals to $(p-1)(1-q)$;
            \item other cases: $(p-1)$ $\mathbb{F}_{q^{p-1}}$-points over $(1:\lambda)$, trace equals to 0.
        \end{itemize}
        Since there are exactly $\frac{q-1}{p-1}-1$ elements $\lambda\in\Fq^\times$ satisfying $(1+\lambda\gamma)^{\frac{q-1}{p-1}}=1$, we find that $\sum_{\lambda\in\Fq^\times}\text{Tr}(\text{Fr}_q; R\Psi_{(1:\lambda)})=(p-1)(q-1)$.

        \item[$(ii)$] The fiber of $(1:0)$ under the resolution of Remark \ref{different resolution} consists of $p-1$ points defined over $\Fq$, where the trace of Frobenius is $1-q$: $\text{Tr}(\text{Fr}_q;R\Psi_{(1:0)})=(p-1)(1-q)$.

        \item[$(iii)$] We look at the fiber of the point $(0:1)$ along $\sigma_2$: we can restrict to $U''_{1,[\tilde{b}=1]}$ and one checks that the $\sigma_2^{-1}((0:1))$ consists of a single point where the trace of Frobenius is $(1-q)^2$ (use that the fiber does not intersect the $[\tilde{v}_0=1]$-chart and that the $[\tilde{v}_1=1]$-chart is semistable if $c$ is invertible). Therefore $\text{Tr}(\text{Fr}_q; R\Psi_{(0:1)})=(1-q)^2$.
        \end{itemize}

        We conclude that $\text{Tr}(\text{Fr}_q; R\Psi_{\gamma})=(1-q)^2$, so the trace is constant on $S_{1,s_1\tau}$.
\end{itemize}

\subsubsection*{Length 0}

This is the most interesting and intricate case. Over the worst singular point of $U$ there is a unique $\Fp$ point given by putting all variables of $U_1'$ to 0: let us denote it by $\uptau_1$. Recall the first three blow-ups in the resolution procedure of Subsection \ref{resolution pro p}:
\[
U_1'''\xrightarrow[]{\sigma_3}U_1''\xrightarrow[]{\sigma_2}U_1'\xrightarrow[]{\sigma_1}U_1.
\]
We need to compute the trace of Frobenius on each point of the fiber $\sigma_1^{-1}(\uptau_1)$, which can be described as in the previous case with the caveat that now $\gamma=0$. We subdivide this fiber in the three cases as before and we find that:
\begin{itemize}
    \item[$(i)$] the fiber over $(1:\lambda)$ the resolution of Remark \ref{different resolution} always consists of $p-1$ $\Fq$-points where the trace is $1-q$, so $\text{Tr}(\text{Fr}_q; R\Psi_{(1:\lambda)})=(p-1)(1-q)$;
    \item[$(ii)$] here the same conclusion as before holds: $\text{Tr}(\text{Fr}_q; R\Psi_{(1:0)})=(p-1)(1-q)$;
    \item[$(iii)$] now the variable $c$ is no longer invertible, and we will need to consider the fibers of the whole chain of blow-ups that form the resolution of $U_1$.
\end{itemize}
We start with $\sigma_2:U_1''\longrightarrow U_1'$. The fiber of $(0:1)$ along this map takes the form of
\[
    \sigma_2^{-1}((0:1))=\text{Proj}(\Fq[\tilde{v}_0,\tilde{v}_1])\begin{cases}
    (1:\delta)\in U''_{1,[\tilde{b}=\tilde{v}_0=1]} \left( \substack{a=b=\tilde{x}=u_1=0\\ v_0=0, \tilde{v}_1=\delta} \right) , \; \delta\in\Fq \\
    (0:1)\in U''_{1,[\tilde{b}=\tilde{v}_1=1]}  \left( \substack{a=b=c=0\\ u_0=v_1=\tilde{v}_0=0} \right).
\end{cases}
\]
There is an isomorphism between $U''_{1,[\tilde{b}=\tilde{v}_1=1]}$ and $U''_{1,[\tilde{b}=\tilde{v}_0=1]}$ sending the point $(0:1)$ to $(1:0)$, so we can restrict our attention to $U''_{1,[\tilde{b}=\tilde{v}_0=1]}$. Recall that $\sigma_3^{-1}(U''_{1,[\tilde{b}=\tilde{v}_0=1]})$ is covered by two charts both isomorphic to the spectrum of
\[
E_0=\frac{\Zp[r,s,t,e_0,f]}{\left( \scalebox{1.2}{$\substack{(r^2st)^{p-1}-p \\ r^{p-1}-e_0f}$} \right)}:
\]
one can check that $\sigma_3^{-1}(1:\delta)$ is a copy of $\mathbb{P}^1_{\Fq}$ and that by intersecting it with one of the two charts we obtain the points $\{e_0=f=r=0, t=\delta, s=\alpha\}_{\alpha\in\Fq}$ inside of $\text{Spec}(E_0)(\Fq)$; call these points $(\alpha,\delta)$. The missing point lies in the other chart, and under its isomorphism with $E_0$ it corresponds to $(0,\delta)$. We are therefore left with the computation of $\text{Tr}(\text{Fr}_q;R\Psi_{(\alpha,\delta)})$ for the nearby cycles on $\text{Spec}(E_0)$: we are going to use the iterated blow-ups described in Step 4 of the previous section. Recall that at the first step we have
\[
\text{Bl}_{(r,e_{0})}(\text{Spec}(E_{0}))=\text{Spec}(E_1)\cup \text{Spec}(R_{0})
\]
with
\[
E_1=\frac{\Zp[r,s,t,e_1,f]}{\left( \scalebox{1.2}{$\substack{(r^2st)^{p-1}-p\\ r^{p-2}-e_1f}$} \right)} \quad \text{and} \quad R_{0}=\frac{\Zp[\tilde{r},s,t,e_{0},f]}{\big( {(\tilde{r}^2e_{0}^2st)^{p-1}-p} \big)}.
\]
Let us denote by $\psi_1:\text{Bl}_{(r,e_{0})}(\text{Spec}(E_{0}))\longrightarrow \text{Spec}(E_{0})$ the projection: the fiber over $(\alpha,\delta)$ is
\[
\psi_1^{-1}((\alpha,\delta))=\text{Proj}(\Fq[\tilde{r},e_1])=\begin{cases}
    (\beta_1:1) \in \text{Spec}(R_0) \left( \substack{e_0=0,\, \tilde{r}=\beta_1\\ s=\alpha,\, t=\delta}\right), \beta_1\in\Fq \\
    (1:0) \in \text{Spec}(E_1) \left( \substack{r=f=e_1=0\\ s=\alpha,\, t=\delta} \right).
\end{cases}
\]
On $\text{Spec}(R_0)$ the trace of Frobenius can be computed directly depending on whether $\alpha,\delta, \beta_1$ are zero or not, while for the trace on $(1:0)$ inside $\text{Spec}(E_1)$ we need to perform another blow-up. We proceed inductively until we reach
\[
\psi_{p-2}:\text{Bl}_{(r,e_{p-3})}(\text{Spec}(E_{p-3}))=\text{Spec}(E_{p-2}) \cup \text{Spec}(R_{p-3})\rightarrow\text{Spec}(E_{p-3}):
\]
at each step the fiber over $(1:0)$ can be described as
\[
\psi_i^{-1}((1:0))=\text{Proj}(\Fq[\tilde{r}, e_i])=\begin{cases}
    (\beta_i:1) \in \text{Spec}(R_{i-1}) \left( \substack{e_i=0,\, \tilde{r}=\beta_i\\ s=\alpha,\, t=\delta} \right), \beta_i\in \Fq \\
    (1:0) \in \text{Spec}(E_i) \left( \substack{r=f=e_i=0\\ s=\alpha,\, t=\delta} \right),
\end{cases}
\]
with the trace of Frobenius directly computable on $R_i$ and on $E_{p-2}$ as well. Summing up all the contributions from the various blow-ups, we find that

\begin{equation}\label{partial trace}
\text{Tr}(\text{Fr}_q, R\Psi_{(\alpha,\delta)}^{E_0})=\sum_{j=1}^{p-2}\Big(\sum_{\beta_j\in\Fq}\text{Tr}(\text{Fr}_q, R\Psi_{(\beta_j:1)}^{R_{j-1}})\Big) + \text{Tr}(\text{Fr}_q; R\Psi_{(1:0)}^{E_{p-2}}).
\end{equation}

By explicit computations, one can check that 
\[
\tr(\fr_q; R\Psi_{(0:1)}^{R_{j-1}})=(1-q)\tr(\fr_q; R\Psi_{(\beta_j:1)}^{R_{j-1}})
\]
for any $\beta_j\in \Fq^\times$. Since $\Fq^\times$ has $q-1$ elements, we see that for $j=1,\ldots,p-2$ we have
\[
\sum_{\beta_j\in\Fq}\text{Tr}(\text{Fr}_q, R\Psi_{(\beta_j:1)}^{R_{j-1}})=0.
\]
Moreover,
\[
E_{p-2}\cong \frac{\Zp[s,t,f,e_{p-2}]}{\big((e_{p-2}^2f^2st)^{p-1}-p\big)},
\]
with the point $(1:0)$ given by $f=e_{p-2}=0, s=\alpha,t=\delta$. It is then easy to check that the expression (\ref{partial trace}) takes the following values:
\begin{itemize}
    \item $\alpha=0, \delta=0$: \quad $\text{Tr}(\text{Fr}_q, R\Psi_{(0,0)}^{E_0})=(1-q)^3$,

    \item $\alpha\neq 0, \delta=0$: \quad $\text{Tr}(\text{Fr}_q, R\Psi_{(\alpha,0)}^{E_0})=(1-q)^2$,

    \item $\alpha=0, \delta\neq 0$: \quad $\text{Tr}(\text{Fr}_q, R\Psi_{(0,\delta)}^{E_0})=(1-q)^2$,

    \item $\alpha\neq 0, \delta\neq 0$: \quad $\text{Tr}(\text{Fr}_q, R\Psi_{(\alpha,\delta)}^{E_0})=(1-q)$.
\end{itemize}

We can now sum everything up and compute the trace of Frobenius over the point $(0:1)$ of $U_1'$. For ease of notation, we will write $\tr(-)$ instead of $\tr(\text{Fr}_q;-)$.
\begin{align*}
&\text{Tr}(R\Psi_{(0:1)}^{U_1'})=2\cdot\text{Tr}(R\Psi_{(1:0)}^{U_1''})+\sum_{\delta\in\Fq^\times}\text{Tr}(R\Psi_{(1:\delta)}^{U_1''})=\\
&=2\Big( 2\cdot\text{Tr}(R\Psi_{(0,0)}^{E_0})+\sum_{\alpha\in\Fq^\times}\text{Tr}(R\Psi_{(\alpha,0)}^{E_0}) \Big) + \sum_{\delta\in\Fq^\times}\Big( 2\cdot\text{Tr}(R\Psi_{(0,\delta)}^{E_0})+ \sum_{\alpha\in\Fq^\times}\text{Tr}(R\Psi_{(\alpha,\delta)}^{E_0}) \Big)=\\
&= 4\cdot\text{Tr}(R\Psi_{(0,0)}^{E_0})+2\sum_{\alpha\in\Fq^\times} \text{Tr}(R\Psi_{(\alpha,0)}^{E_0}) +2\sum_{\delta\in\Fq^\times} \text{Tr}(R\Psi_{(0,\delta)}^{E_0}) + \sum_{\alpha,\delta\in\Fq^\times} \text{Tr}(R\Psi_{(\alpha,\delta)}^{E_0})=\\
&=(1-q)^3.
\end{align*}

Finally, the trace of Frobenius on the preimage of the worst singular point is:
\begin{align*}
    &\text{Tr}(\text{Fr}_q, R\Psi_{\tau_1})=\sum_{\lambda\in\Fq^\times}\text{Tr}(\text{Fr}_q, R\Psi_{(1:\lambda)}^{U_1'}) + \text{Tr}(\text{Fr}_q, R\Psi_{(1:0)}^{U_1'}) + \text{Tr}(\text{Fr}_q, R\Psi_{(0:1)}^{U_1'})=\\
    &=(1-q)^3+(p-1)q(1-q).
\end{align*}

\section{Comparison with the test function}\label{Comparison with the test function}

In this section, we compare the semisimple trace of Frobenius on the pushforward $\pi_*R\Psi^1$ with the values of Horn's test function of \cite[Example 5.4.1]{Hor}, in the case of $\text{GSp}_4$. Before doing this, we recall the work of Haines and Rapoport on some unitary Shimura varieties of the kind considered in this thesis, since the situation in the Siegel case presents some similarities but also different aspects.

\subsection{Work of Haines-Rapoport on the ``good'' Drinfeld case}\label{HR drinfeld case}

In \cite{HR}, Haines and Rapoport prove that the conjectural test function proposed by Haines and Kottwitz in \cite{Hai} is actually a test function in the case of unitary Shimura varieties considered in this thesis, when the signature is $(n-1,1)$ (the ``good'' Drinfeld case). In this case both the integral models $\mathcal{A}_0$ and $\mathcal{A}_1$ are semistable and the picture on the local model side is 
\[
U_1=\text{Spec}\left( \frac{\Zp[u_0,\ldots,u_{n-1}]}{((u_0\cdots u_{n-1})^{p-1}-p} \right) \xrightarrow{\makebox[1.3cm]{\tiny{$u_i^{p-1}\mapsfrom x_i$}}} \text{Spec}\left( \frac{\Zp[x_0,\ldots,x_{n-1}]}{(x_0\cdots x_{n-1}-p)} \right)=U.
\]
They prove that the semi-simple Lefschetz number
\[
\text{Lef}^{ss}(\text{Fr}_q; R\Psi^1)=\sum_{x\in\mathcal{A}_1(\Fq)}\tr^{ss}(\fr_q;R\Psi^1_x)
\]
can be written as a certain sum of the form
\[
\sum_{(\gamma_0;\gamma,\delta)}c(\gamma_0;\gamma,\delta)\text{O}_\gamma(f^p)\text{TO}_{\delta\sigma}(\phi_{r,1}),
\]
where $\phi_{r,1}$ is the test function of Kottwitz and Haines in this specific pro-$p$ Iwahori case. They use the equality $\text{Lef}^{ss}(\text{Fr}_q; R\Psi^1)=\text{Lef}^{ss}(\fr_q;\pi_*R\Psi^1)$, where $\pi:\mathcal{A}_1\longrightarrow\mathcal{A}_0$ is the projection, thanks to which it is enough to consider the push-forward of the nearby cycles complex. Moreover, $\pi$ is finite flat of degree $(p-1)^n$ and it is a connected Galois cover on the generic fiber with Galois group $T(\Fp)$, $T$ being the diagonal maximal torus of $\text{GL}_n$, so that
\[
\pi_{\eta*}(\overline{\Q}_\ell)=\bigoplus_{T(\Fp)\xrightarrow[]{\chi}\Qb_\ell^\times} \Qb_{\ell,\chi}^\times,
\]
where $\Qb_{\ell,\chi}^\times$ is the rank 1 local system on $\mathcal{A}_{0,\eta}$ corresponding to $\pi_1(\mathcal{A}_{0,\eta})\twoheadrightarrow T(\Fp)\xrightarrow[]{\chi}\Qb_\ell^\times$. Setting $R\Psi_\chi=R\Psi(\Qb_{\ell,\chi}^\times)$, by compatibility of nearby cycles with proper maps they obtain a decomposition
\[
\pi_*R\Psi^1=\bigoplus_{T(\Fp)\xrightarrow[]{\chi}\Qb_\ell^\times}R\Psi_\chi
\]
and they prove the following.

\begin{theorem}[\cite{HR}]
    Let $x\in\mathcal{A}_{0,w}(\Fq)$. Then for any element $t\in T(\Fq)$ such that $N_r(t)$ projects to $t_x\in T^{S(w)}(\Fp)$, we have
    \[
    \emph{Tr}(\emph{Fr}_q;(R\Psi_\chi)_x)=\phi_{r,\chi}(t^{-1}w^{-1}).
    \]
\end{theorem}
Let us explain the notation used in this theorem: $\phi_{r,\chi}$ is a specific element inside a suitable pro-$p$ Iwahori-Hecke algebra (see \cite[Definition 7.1.2]{HR}), $N_r:\Fq^\times \rightarrow \Fp^\times$ is the norm and $w\in \text{Adm}(\mu)$ is a $\mu$-admissible element for the cocharacter $\mu=(1^{(n-1)},0)$ of $\text{GL}_n$ with associate KR-stratum $\mathcal{A}_{0,w}$. Note that the KR-stratification can be indexed by subsets $S\subset\{1,\ldots,n\}$ ($S$ corresponds to the points of $\mathcal{A}_{0,\Fp}$ where the groups $G_i$ are infinitesimal for $i\in S$, see \cite[Section 3.4]{HR}) and $S(w)$ is the subset corresponding to the stratum associated to $w$. Given such a subset $S$ one can consider the short exact sequence of tori
\[
1\longrightarrow T_S\longrightarrow T\longrightarrow T^S\longrightarrow1
\]
where $T_S=\prod_{i\in S}\mathbb{G}_m$ and $T^S=\prod_{i\notin S}\mathbb{G}_m$: it turns out that the restriction of $\pi$ to $(\mathcal{A}_{1,w})_\text{red}\longrightarrow\mathcal{A}_{0,w}$ is a $T^{S(w)}$-torsor, see \cite[Corollary 3.4.4]{HR}; here $\mathcal{A}_{1,w}=\pi^{-1}(\mathcal{A}_{0,w})$. Now, for $x\in \mathcal{A}_{0,w}(\Fq)$ we can choose $y\in\pi^{-1}(x)$: since $\fr_q$ permutes the fiber over $x$, there exists a unique element $t_x\in T^{S(w)}$ defined by the equation $\fr_q(y)=t_x\cdot y$ and this element does not depend on the choice of $y$ (see \cite[Section 7.3]{HR}). Haines and Rapoport then define 
\[
\phi_{r,1}=(q-1)^{-n}\sum_{\chi\in T(\Fp)^\vee}\phi_{r,\chi},
\]
prove that it lies in the center of the pro-$p$ Iwahori-Hecke algebra and compute it explicitly in \cite[Proposition 12.2.1]{HR}, which gives
\begin{equation}\label{pro-p Iw test function Drinfeld}
\phi_{r,1}(tw^{-1})=\begin{cases}
    0, \quad \text{if}\; w\notin\text{Adm}(\mu),\\
    0, \quad \text{if}\; w\in\text{Adm}(\mu)\; \text{but}\; N_r(t)\notin T_{S(w)}(\Fp),\\
    (-1)^n(p-1)^{n-|S(w)|}(1-q)^{|S(w)|-n-1}, \quad \text{otherwise}.
\end{cases}
\end{equation}

\begin{remark}
One can also show directly that $\phi_{r,1}$ computes the trace of Frobenius on $\pi_*R\Psi^1$ (up to the factor $(q-1)^{-n}$): for $x\in \mathcal{A}_{0,w}(\Fq)$ we have
\begin{equation}\label{stalk HR}
(\pi_*R\Psi^1)_x=\bigoplus_{y\overset{\pi}{\mapsto}x}R\Psi_y^1.
\end{equation}
Given such an $x$, fix a trivialization of the Oort-Tate line bundles $\mathcal{L}_i$ of the groups $G_i$: we can see the Oort-Tate parameters of the groups $G_i$ as pairs $(a_i,b_i)\in(\Fq)^2$ such that $a_ib_i=0$, and the fiber of $\pi$ over $x$ is given by extracting $(p-1)$-th roots $u_i\in\overline{\mathbb{F}}_p$ of each $a_i$. Note that $a_i=0$ if and only if $i\in S(w)$ and that if $i\notin S(w)$ any two roots of $a_i$ differ by an element in $\Fp^\times$. Therefore, since $\fr_q(y)=y$ if and only if $y\in \mathcal{A}_1(\Fq)$ for $y\in \pi^{-1}(y)$, there are two disjoint cases:
\medskip
\begin{itemize}
    \item[$(i)$] $\pi^{-1}(x)\cap \mathcal{A}_1(\Fq)=\emptyset$, if and only if $t_x\neq 1$: $\fr_q$ permutes the summands of (\ref{stalk HR}) and its trace on $(\pi_*R\Psi^1)_x$ is 0;
    \item[$(ii)$] $\pi^{-1}(x)\subset \mathcal{A}_1(\Fq)$, if and only if $t_x=1$: $\fr_q$ acts diagonally on $(\pi_*R\Psi^1)_x$ and its trace on $R\Psi_y^1$ is equal to $(1-q)^{|S(w)|-1}$.
\end{itemize}
Moreover, the cardinality of $\pi^{-1}(x)$ is always $(p-1)^{n-|S(w)|}$, so we find that
\[
\tr(\fr_q; (\pi_*R\Psi)_x)=\begin{cases}
    0, \quad \text{if} \; t_x\neq 1;\\
    (p-1)^{n-|S(w)|}(1-q)^{|S(w)|-1} \quad \text{if} \; t_x=1.
\end{cases}
\]
Here we are using the fact that $\ApI$ is semistable, with special fiber a non-reduced divisor whose components have multiplicity $p-1$, so that the trace of $\fr$ at any $\Fq$-point of $\mathcal{A}_{1,w}$ is $(1-q)^{|S(w)|-1}$, as in the Iwahori case. Now, given $t\in T(\Fq)$ such that $N_r(t)$ projects to $t_x$, we have that $N_r(t)\in T_{S(w)}(\Fp)$ if and only if $t_x=1$, because $T_{S(w)}$ is the kernel of the projection $T\rightarrow T^{S(w)}$. Therefore we find that
\[
\tr(\fr_q;(\pi_*R\Psi^1)_x)= (q-1)^n\phi_{r,1}(tw^{-1}),
\]
where $w\in\text{Adm}(\mu)$ corresponds to the KR-stratum of $x$ and $t$ is as above. However, with this approach it is not clear that this function $\phi_{r,1}$ lies in the center of the pro-$p$ Iwahori-Hecke algebra: to get this result, it is crucial to decompose it into its $\chi$-components $\phi_{r,\chi}$ and proceed as Haines and Rapoport do.
\end{remark}

\begin{remark}\label{specific element}
    We can find an element $s_x\in T(\Fq)$, which depends only on the point $x\in\mathcal{A}_{0,w}(\Fq)$ and on a trivialization of the associated Oort-Tate line bundle $\mathcal{L}_i$, such that 
    \begin{equation}\label{test function comp Drinfeld}
    \tr(\fr_q;(\pi_*R\Psi^1)_x)= (q-1)^n\phi_{r,1}(s_xw^{-1}).
    \end{equation}
    Consider the Oort-Tate parameters $(a_0,\ldots,a_{n-1})$ associated with $x$: then we set
    \[
    s_x=\scalebox{0.8}{$
    \begin{pmatrix}
        a_0' & &\\
        & \ddots &\\
        & & a_{n-1}'
    \end{pmatrix}$ }, \quad \text{where} \;\; 
    a_i'=\begin{cases}
        a_i,  &\text{if $a_i\neq 0$};\\
        1,  &\text{if $a_i=0$}.
    \end{cases}
    \]
    Now consider $N_r(s_x)\in T(\Fp)$: we claim that it projects to $t_x$ under $T\rightarrow T^{S(w)}$. Indeed, $t_x=\text{diag}(t_x^{i})_{i\notin S(w)}$ and by its definition it has to satisfy the relation $t_x^iu_i=u_i^q$ for $i\notin S(w)$, where $y=(x; u_i)_{i=0}^{n-1}$ is a point of $\mathcal{A}_{1,w}$ lying over $x$. Thus we have that 
    \[
    t_x^i=u_i^{q-1}=a_i^{\frac{q-1}{p-1}}=N_r(a_i),
    \]
    the fact that $a_i\neq 0$ if and only if $i\neq S(w)$ proves the claim and the formula (\ref{test function comp Drinfeld}) follows.
\end{remark}

To summarize, to any $ w \in\text{Adm}(\mu)$ one can associate the subgroup $T_{S(w)}\subset T(\Fp)$ and then consider the function in the pro-$p$ Iwahori-Hecke algebra defined by the formula (\ref{pro-p Iw test function Drinfeld}). Now, for any $x\in\mathcal{A}_{0,w}(\Fq)$ we can find an element $s_x\in T(\Fq)$ and use it to compare the test function $\phi_{r,1}$ with the trace of $\fr_q$ on $\pi_*R\Psi^1$. In the $\text{GSp}_4$ case the situation is somewhat similar, but with some differences that we are going to highlight: most notably, the $p$-rank zero locus will not be a single point and the procedure for finding $s_x$ will not always give an element for which a formula like (\ref{test function comp Drinfeld}) holds.

\subsection{The \texorpdfstring{$\text{GSp}_4$}{} case}\label{the GSp4 case}

For this subsection, we denote by $\AI$ the integral model for $\text{GSp}_4$ in the Iwahori case and by $\ApI$ the Haines-Stroh integral model for $\text{GSp}_4$ in the pro-$p$ Iwahori case. We are going to compare the trace of Frobenius on the pushforward of $R\Psi^1$ along $\pi:\ApI\longrightarrow \AI$ with the test function $\phi_{r,1}'$ of \cite[Example 5.4.1]{Hor}. It is a function in the pro-$p$ Iwahori-Hecke algebra for $\text{GSp}_4$ so it can be described by its values on elements $sw$ of the pro-$p$ extended Weyl group $\widetilde{W}_1=T(\Fq)\rtimes \widetilde{W}$. If $w\notin\text{Adm}(\mu)$, then $\phi_{r,1}'(sw)=0$ for any $s\in T(\Fq)$; if instead $w\in\text{Adm}(\mu)$, the value of $\phi_{r,1}'$ on $sw$ can be either $0,1$ or a polynomial in $p$ and $q$ depending on a condition on $s$, namely if $s$ belongs to a certain subgroup $A_{w,J_{\Delta_1}}\subset T(\Fp)$, defined in \cite[Definition 3.27]{Hor} (actually for $w=\tau$ the situation is slightly different, see below). We collect these subgroups in the table below. In the middle column, by $\text{diag}(\alpha,\beta,\alpha,\beta)$ we mean the subgroup of $T(\Fp)$ consisting of diagonal matrices of this form, with $\alpha,\beta\in\Fp^\times$ (and similarly for the other cases).

\begin{center}
\renewcommand{\arraystretch}{1.25}
\begin{tabular}{|p{3cm}|p{3cm}|p{2.8cm}|}
\hline
$\text{Adm}(\mu)$ & Subgroup $A_{w,J_{\Delta_1}}$ & OT parameters \\
\hline
$s_{010}\tau=t_{(1100)}$ & $\text{diag}(\alpha,\alpha,1,1)$ &$(0,0,a_1,a_0)$ \\ 
$s_{102}\tau=t_{(0101)}$ & $\text{diag}(1,\alpha,1,\alpha)$  &$(b_0,0,a_1,0)$\\
$s_{201}\tau=t_{(1010)}$ & $\text{diag}(\alpha,1,\alpha,1)$ &$(0,b_1,0,a_0)$ \\
$s_{212}\tau=t_{(0011)}$ & $\text{diag}(1,1,\alpha,\alpha)$ &$(b_0,b_1,0,0)$ \\
$s_{01}\tau=t_{(1100)}s_{2}$ & $\text{diag}(\alpha\beta,\alpha,\beta,1)$ &$(0,0,0,a_0)$ \\
$s_{12}\tau=t_{(0101)}s_{2}$ & $\text{diag}(1,\alpha,\beta,\alpha\beta)$ &$(b_0,0,0,0)$ \\
$s_{10}\tau=t_{(1100)}s_{121}$ & $\text{diag}(\alpha,\alpha\beta,1,\beta)$ &$(0,0,a_1,0)$ \\
$s_{02}\tau=t_{(1010)}s_{1}$ & $\text{diag}(\alpha,\beta,\alpha,\beta)$ &$(0,0,0,0)$ \\
$s_{21}\tau=t_{(1010)}s_{121}$ & $\text{diag}(\alpha,1,\alpha\delta,\delta)$ &$(0,b_1,0,0)$ \\
$s_{0}\tau=t_{(1100)}s_{21}$ & $T(\Fp)$ &$(0,0,0,0)$ \\
$s_{1}\tau=t_{(1100)}s_{1212}$ & $T(\Fp)$ &$(0,0,0,0)$ \\
$s_{2}\tau=t_{(1010)}s_{12}$ & $T(\Fp)$ &$(0,0,0,0)$ \\
$\tau=t_{(1100)}s_{212}$ & $\text{diag}(\alpha,\alpha,\beta,\beta)$ &$(0,0,0,0)$ \\
\hline
\end{tabular}
\end{center} 

They are the analogue of the subtori $T^S(\Fp)$ and in most cases they are reflected by the geometry of the morphism $\pi$, the only exceptions being the worst singular point and the element $s_{02}\tau$. Note that in Horn's thesis there is a second subgroup appearing for $\tau$, $A_{\tau,J_{\Delta_2}}$: however this group is trivial, so we can ignore it. We now explain the relation of the subgroups $A_{w,J_{\Delta_1}}$ with the geometry of the pro-$p$ integral model. Recall that for $x\in\AI$ we have the group schemes $G_i=\text{Ker}(\alpha_i)$ with Oort-Tate parameters $(\mathcal{L}_i,a_i,b_i)$, and points of $\ApI$ are tuples $(x; v_0,v_1,u_1,u_0)$ where $x\in\AI$ and $v_i,u_i$ are Oort-Tate generators for $G^*_i$ and $G_i$ respectively, such that $u_0v_0=u_1v_1$. The $\Fp$-points of the diagonal torus of $\text{GSp}_4$
\[
T= \{ \underline{t}= \text{diag}(h_0,h_1,k_1,k_0) \,|\, h_0k_0=h_1k_1 \}
\]
act on $\ApI$ by 
\[
\underline{t}\cdot(x; v_0,v_1,u_1,u_0)=(x; h_0v_0,h_1v_1,k_1u_1,k_0u_0)
\]
and $\pi$ is invariant under this action. Given $x\in\mathcal{A}_{0,w}(\overline{\mathbb{F}}_p)$, let $e_x=(b_0,b_1,a_1,a_0)$ be the vector of its Oort-Tate parameters: its 0 entries do not depend on $x$, but only on $w$. We can consider the short exact sequence of tori
\[
1\rightarrow T_w\longrightarrow T \longrightarrow T^w \rightarrow 1,
\]
where $T_w$ is the subtorus of $T$ consisting of those $\underline{t}\in T$ having 1 in the entries corresponding to the non-0 coordinates of $e_x$, for some $x\in\mathcal{A}_{0,w}(\overline{\mathbb{F}}_p)$. By construction, $T_w(\Fp)$ acts trivially on $\mathcal{A}_{1,w}$ and as in the case considered by Haines and Rapoport we see that the restriction of $\pi$ to $(\mathcal{A}_{1,w})_\text{red}\rightarrow \mathcal{A}_{0,w}$ is a torsor for the group $T^w(\Fp)$. From the table above it is clear that the two subgroups $T_w$ and $A_{w,J_{\Delta_1}}$ coincide in all cases except $\tau$ and $s_{02}\tau$: we will see that on the worst singular point this does not cause problems, whereas for $s_{02}\tau$ the situation is more complicated. Now, for each $x\in \mathcal{A}_{0,w}(\Fq)$ we can define an element $\tilde{s}_x\in T(\Fq)$ exactly as in Remark \ref{specific element}: it turns out that this element $\tilde{s}_x$, as in the case considered by Haines and Rapoport, can be used to compare Horn's test function with the trace of Frobenius on nearby cycles except for the stratum corresponding to $s_{02}\tau$; in this case, a more ad hoc definition is needed. We can also define $t_x\in T^w(\Fp)$ following the procedure of Haines and Rapoport and check that $N_r(\tilde{s}_x)$ projects to $t_x$ under $T\rightarrow T^w$. We now compare the test function with the trace of Frobenius on nearby cycles.

\subsubsection*{Length 3}

In this case, for $w=t_\lambda\in\text{Adm}(\mu)$ of length 3, Horn's test function takes the following form:
\[
\phi'_{r,1}(sw)=\begin{cases}
    0, & \text{if}\; N_r(s)\notin A_{w,J_{\Delta_1}};\\
    (-1)(p-1)^2(1-q)^{-3}, & \text{if}\; N_r(s)\in A_{w,J_{\Delta_1}}.
\end{cases}
\]
On the other hand, we have that for $x\in\mathcal{A}_{0,w}(\Fq)$ 
\[
\tr(\fr_q; (\pi_*R\Psi^1)_x)=\begin{cases}
    0, &\text{if} \; t_x\neq 1;\\
    (p-1)^2, &\text{if}\; t_x=1.
\end{cases}
\]
Following the same reasoning of Remark \ref{specific element} we obtain that
\[
\tr(\fr_q; (\pi_*R\Psi^1)_x)=(q-1)^3\phi'_{r,1}(\tilde{s}_xw).
\]

\subsubsection*{Length 2}

For $w\in\{s_{12}\tau, s_{10}\tau, s_{21}\tau, s_{01}\tau\}\subset\text{Adm}(\mu)$, Horn's test function takes the following form:
\[
\phi'_{r,1}(sw)=\begin{cases}
    0, & \text{if}\; N_r(s)\notin A_{w,J_{\Delta_1}};\\
    (-1)(p-1)(1-q)^{-2}, & \text{if}\; N_r(s)\in A_{w,J_{\Delta_1}}.
\end{cases}
\]
On the geometric side, for $x\in\mathcal{A}_{0,w}(\Fq)$ we have
\[
\tr(\fr_q; (\pi_*R\Psi^1)_x)=\begin{cases}
    0, &\text{if} \; t_x\neq 1;\\
    (p-1)(1-q), &\text{if}\; t_x=1.
\end{cases}
\]
So we obtain again that
\[
\tr(\fr_q; (\pi_*R\Psi^1)_x)=(q-1)^3\phi'_{r,1}(\tilde{s}_xw).
\]

If $w=s_{02}\tau$, Horn's formula is still of the same form, but in this case the group $A_{w,J_{\Delta_1}}$ is different from $T^w(\Fp)$, which is trivial: this is because the stratum $S_{s_{02}\tau}$ is contained in the $p$-rank zero locus, over which the groups $T^w$ are always trivial. Nevertheless, from the formula (\ref{trace s02}) it follows that 
\[
\tr(\fr_q; (\pi_*R\Psi^1)_x)=\begin{cases}
    0, &\text{if} \; N_r(\alpha/\beta)\neq 1;\\
    (p-1)(1-q), &\text{if}\; N_r(\alpha/\beta)= 1.
\end{cases}
\]
Here $(\alpha, \beta)$ is a point of the open $U\subset \text{M}_1^{\text{GSp}_4}$ corresponding to $x\in\mathcal{A}_{0,w}$ under the local model diagram (see also the discussion in Subsection \ref{trace of frob on U1} relative to $s_{02}\tau$). Of course such a point is not unique, but the trace of $\fr_q$ on $R\Psi_{(\alpha,\beta)}^{U}$ only depends on $x$. Now, consider the element $s_x=\text{diag}(\alpha,\beta,\alpha,\beta)\in T(\Fq)$: it holds that $N_r(s_x)\in A_{s_{02}\tau, J_{\Delta_1}}$ if and only if $N_r(\alpha/\beta)=1$. Therefore we obtain again
\[
\tr(\fr_q; (\pi_*R\Psi^1)_x)=(q-1)^3\phi'_{r,1}(s_xw).
\]
for any choice of $s_x$ as above. It is slightly unpleasant that this element $s_x$ does not depend only on $x$, but also on the choice of corresponding point in the local model. However, we could not find a way of characterizing such an $s_x$ directly from $x$.

\subsubsection*{Length 1}

For $w\in\text{Adm}(\mu)$ of length 1, Horn's test function is
\[
\phi'_{r,1}(sw)= (-1)(1-q)^{-1}
\]
for any $s\in T(\Fq)$, since $A_{w,J_{\Delta_1}}$ is equal to $T(\Fq)$ (note that for $w=t_{(1100)}s_{2121}=s_1\tau$ there is a small mistake and $J_{\Delta_1}$ should be $\{ \pm \alpha_2, \pm\alpha_4 \}$). For $x\in\mathcal{A}_{0,w}(\Fq)$ we find that 
\[
\tr(\fr_q; (\pi_*R\Psi^1)_x)=(1-q)^2
\]
and therefore
\[
\tr(\fr_q; (\pi_*R\Psi^1)_x)=(q-1)^3\phi'_{r,1}(\tilde{s}_xw).
\]

\subsubsection*{Length 0}

For $w=\tau$ Horn's test function is
\[
\phi_{r,1}'=\begin{cases}
    -1, & \text{if $N_r(s)\notin A_{w, J_{\Delta_1}}$}\\
    (-1)\big( 1+(p-1)q(1-q)^{-2} \big) & \text{if $N_r(s)\in A_{w, J_{\Delta_1}}$}
\end{cases}
\]
and the trace of Frobenius is
\[
\tr(\fr_q; (\pi_*R\Psi^1)_x))= (1-q)^3+(p-1)q(1-q),
\]
thus $\tr(\fr_q; (\pi_*R\Psi^1)_x)= (q-1)^3\phi_{r,1}'(\tilde{s}_x\tau)$ (since in this case $s_x=1$).

\bigskip

\begin{notation}
    Let $x\in\mathcal{A}_{0,w}(\Fq)$: if $w\neq s_{02}\tau$, we write $s_x\in T(\Fq)$ for the element $\tilde{s}_x$; if $w=s_{02}\tau$, we write $s_x$ for the element defined above in this case.
\end{notation}

In conclusion, we have proved the following theorem.

\begin{theorem}\label{comparison test function}
    For all $x\in \mathcal{A}_0(\Fq)$ we have that
    \[
    \emph{Tr}(\emph{Fr}_q; (\pi_*R\Psi^1)_x)=(q-1)^3\phi'_{r,1}(s_xw),
    \]
    where $w\in \emph{Adm}(\mu)$ corresponds to the KR-stratum which $x$ belongs to. 
\end{theorem}

\end{document}